\documentclass[12pt,a4paper]{article}
\usepackage{amsmath}
\usepackage{amssymb}
\usepackage{enumerate}
\usepackage{t1enc}
\usepackage[latin1]{inputenc}
\usepackage[english]{babel}
\usepackage{enumerate}

\usepackage{mathtools}

\usepackage{amsfonts}
\usepackage{latexsym}
\usepackage{bm}
\newtheorem{theorem}{Theorem} 

\newtheorem{lemma}{Lemma}

\newtheorem{conj}{Conjecture}

\newtheorem{construction}{Construction}
\newtheorem{problem}{Problem}

\allowdisplaybreaks

\begin{document}
\title{Isolation of connected graphs}

\author{Peter Borg\\[5mm]
{\normalsize Department of Mathematics} \\
{\normalsize Faculty of Science} \\
{\normalsize University of Malta}\\
{\normalsize Malta}\\
{\normalsize \texttt{peter.borg@um.edu.mt}}
}

\date{}
\maketitle

\begin{abstract} 
For a connected $n$-vertex graph $G$ and a set $\mathcal{F}$ of graphs, let $\iota(G,\mathcal{F})$ denote the size of a smallest set $D$ of vertices of $G$ such that the graph obtained from $G$ by deleting the closed neighbourhood of $D$ contains no graph in $\mathcal{F}$. Let $\mathcal{E}_k$ denote the set of connected graphs that have at least $k$ edges. By a result of Caro and Hansberg, $\iota(G,\mathcal{E}_1) \leq n/3$ if $n \neq 2$ and $G$ is not a $5$-cycle. The author recently showed that if $G$ is not a triangle and $\mathcal{C}$ is the set of cycles, then $\iota(G,\mathcal{C}) \leq n/4$. We improve this result by showing that $\iota(G,\mathcal{E}_3) \leq n/4$ if $G$ is neither a triangle nor a $7$-cycle. Let $r$ be the number of vertices of $G$ that have only one neighbour. We determine a set $\mathcal{S}$ of six graphs such that $\iota(G,\mathcal{E}_2) \leq (4n - r)/14$ if $G$ is not a copy of a member of $\mathcal{S}$. The bounds are sharp.
\end{abstract}

\section{Introduction}
Unless stated otherwise, we use small letters such as $x$ to denote non-negative integers or elements of sets, and capital letters such as $X$ to denote sets or graphs. The set of positive integers is denoted by $\mathbb{N}$. For $n \geq 1$, $[n]$ denotes the set $\{1, \dots, n\}$ (that is, $[n]= \{i \in \mathbb{N} \colon i \leq n\}$). We take $[0]$ to be the empty set $\emptyset$. Arbitrary sets are taken to be finite. For a set $X$, ${X \choose 2}$ denotes the set of $2$-element subsets of $X$ (that is, ${X \choose 2} = \{ \{x,y \} \colon x,y \in X, x \neq y \}$). We may represent a $2$-element set $\{x,y\}$ by $xy$.

For standard terminology in graph theory, we refer the reader to \cite{West}. Most of the notation and terminology used here is defined in \cite{Borg1}, which lays the foundation for the work presented here. 

Every graph $G$ is taken to be \emph{simple}, that is, $G$ is a pair $(V(G), E(G))$ such that $V(G)$ and $E(G)$ (the vertex set and the edge set of $G$) are sets that satisfy $E(G) \subseteq {V(G) \choose 2}$. We call $G$ an \emph{$n$-vertex graph} if $|V(G)| = n$. We call $G$ an \emph{$m$-edge graph} if $|E(G)| = m$. For a vertex $v$ of $G$, $N_{G}(v)$ denotes the set $\{w \in V(G) \colon vw \in E(G)\}$ of neighbours of $v$ in $G$, $N_{G}[v]$ denotes the closed neighbourhood $N_{G}(v) \cup \{ v \}$ of $v$, and $d_{G}(v)$ denotes the degree $|N_{G} (v)|$ of $v$. 
For a subset $X$ of $V(G)$, $N_G[X]$ denotes the closed neighbourhood $\bigcup_{v \in X} N_G[v]$ of $X$, $G[X]$ denotes the subgraph of $G$ induced by $X$ (that is, $G[X] = (X,E(G) \cap {X \choose 2})$), and $G - X$ denotes the graph obtained by deleting the vertices in $X$ from $G$ (that is, $G - X = G[V(G) \backslash X]$). Where no confusion arises, the subscript $G$ may be omitted from any notation that uses it; for example, $N_G(v)$ may be abbreviated to $N(v)$. If $H$ is a subgraph of $G$, then we say that \emph{$G$ contains $H$}. 
A \emph{component of $G$} is a maximal connected subgraph of $G$ (that is, no other connected subgraph of $G$ contains it). Clearly, the components of $G$ are pairwise vertex-disjoint (that is, no two have a common vertex), and their union is $G$. If $F$ is a copy of $G$, then we write $F \simeq G$.
  
For $n \geq 1$, the graphs $([n], {[n] \choose 2})$ and $([n], \{\{i,i+1\} \colon i \in [n-1]\})$ are denoted by $K_n$ and $P_n$, respectively. For $n \geq 3$, $C_n$ denotes the graph $([n], \{\{1,2\}, \{2,3\}, \dots, \{n-1,n\}, \{n,1\}\})$ ($= ([n], E(P_n) \cup \{\{n,1\}\})$). A copy of $K_n$ is called an \emph{$n$-clique} or a \emph{complete graph}. A copy of $P_n$ is called an \emph{$n$-path} or simply a \emph{path}. A copy of $C_n$ is called an \emph{$n$-cycle} or simply a \emph{cycle}. We call a $3$-cycle a \emph{triangle}. Note that $K_3$ is the triangle $C_3$. 

If $\mathcal{F}$ is a set of graphs and $F$ is a copy of a graph in $\mathcal{F}$, then we call $F$ an \emph{$\mathcal{F}$-graph}. If $G$ is a graph and $D \subseteq V(G)$ such that $G-N[D]$ contains no $\mathcal{F}$-graph, then $D$ is called an \emph{$\mathcal{F}$-isolating set of $G$}. The size of a smallest $\mathcal{F}$-isolating set of $G$ is denoted by $\iota(G, \mathcal{F})$ and called the \emph{$\mathcal{F}$-isolation number of $G$}. We abbreviate $\iota(G, \{F\})$ to $\iota(G, F)$.  

The study of isolating sets was initiated recently by Caro and Hansberg~\cite{CaHa17}. It generalizes the study of the classical domination problem \cite{C, CH, HHS, HHS2, HL, HL2} naturally. Indeed, $D$ is a \emph{dominating set of $G$} (that is, $N[D] = V(G)$) if and only if $D$ is a $\{K_1\}$-isolating set of $G$, so the \emph{domination number} is the $\{K_1\}$-isolation number. One of the earliest domination results is the upper bound $n/2$ of Ore \cite{Ore} on the domination number of any connected $n$-vertex graph $G \not\simeq K_1$ (see \cite{HHS}). While deleting the closed neighbourhood of a dominating set produces the graph with no vertices, deleting the closed neighbourhood of a $\{K_2\}$-isolating set produces a graph with no edges. A major contribution of Caro and Hansberg in \cite{CaHa17} is a sharp upper bound for this case. They proved that if $G$ is a connected $n$-vertex graph, then $\iota(G, K_2) \leq n/3$ unless $G \simeq K_2$ or $G \simeq C_5$. Fenech, Kaemawichanurat, and the present author \cite{BFK} generalized these bounds by showing that for any $k \geq 1$, $\iota(G, K_k) \leq n/(k+1)$ unless $G \simeq K_k$ or $k=2$ and $G$ is a $5$-cycle. This sharp bound settled a problem of Caro and Hansberg~\cite{CaHa17}.

Let $\mathcal{C}$ denote the set of cycles. The author \cite{Borg1} obtained a sharp upper bound on $\iota(G, \mathcal{C})$, and consequently settled another problem of Caro and Hansberg \cite{CaHa17}. Before stating the result, we recall the explicit construction used for an extremal case. 

\begin{construction}[\cite{Borg1}] \label{Bconstruction} \emph{For any $n, k \in \mathbb{N}$ and any $k$-vertex graph $F$, we construct a connected $n$-vertex graph $B_{n,F}$ as follows.  If $n \leq k$, then let $B_{n,F} = P_n$. If $n \geq k+1$, then let $a_{n,k} = \left \lfloor \frac{n}{k+1} \right \rfloor$, let $b_{n,k} = n - ka_{n,k}$ (so $a_{n,k} \leq b_{n,k} \leq a_{n,k} + k$), let $F_1, \dots, F_{a_{n,k}}$ be copies of $F$ such that $P_{b_{n,k}}, F_1, \dots, F_{a_{n,k}}$ are pairwise vertex-disjoint, and let $B_{n,F}$ be the graph with $V(B_{n,F}) = [b_{n,k}] \cup \bigcup_{i=1}^{a_{n,k}} V(F_i)$ and $E(B_{n,F}) = E(P_{a_{n,k}}) \cup \{\{a_{n,k}, j\} \colon j \in [b_{n,k}] \backslash [a_{n,k}]\} \cup \bigcup_{i=1}^{a_{n,k}} (E(F_i) \cup \{\{i,v\} \colon v \in V(F_i)\})$.}
\end{construction}
\begin{theorem}[\cite{Borg1}] \label{Borgcycle}
If $G$ is a connected $n$-vertex graph that is not a triangle, then
\[\iota(G, \mathcal{C}) \leq \left \lfloor \frac{n}{4} \right \rfloor.\] 
Moreover, equality holds if $G = B_{n,K_3}$.
\end{theorem}

For $k \geq 1$, let $\mathcal{E}_k$ denote the set of connected graphs that have at least $k$ edges. For a graph $G$, we may abbreviate $\iota(G,\mathcal{E}_k)$ to $\iota_k(G)$. For $j \geq k$, $\iota_j(G) \leq \iota_k(G)$ as every $\mathcal{E}_k$-isolating set of $G$ is an $\mathcal{E}_{j}$-isolating set of $G$. 
In this paper, we introduce the problem of determining a sharp upper bound on $\iota_k(G)$ similar to the bounds mentioned above. The above-mentioned result of Caro and Hansberg \cite{CaHa17} solves the problem for $k=1$. We solve the problem for $k=2$ and for $k=3$.

\begin{theorem}[\cite{CaHa17}] \label{CHresult} If $G$ is a connected $n$-vertex graph that is neither a $2$-clique nor a $5$-cycle, then 
\[\iota_1(G) \leq \left \lfloor \frac{n}{3} \right \rfloor.\] 
%
\end{theorem}
The bound in Theorem~\ref{CHresult} is attained if $G = B_{n,K_2}$.

For $k \geq 1$, let $K_{1,k}$ denote the star $([k+1], \{\{1,i\} \colon i \in [k+1] \backslash \{1\}\})$. Let $C_6'$ denote the graph $([7], E(C_6) \cup \{\{1,7\}\})$. Let $C_6''$ denote the graph $([7], E(C_6') \cup \{\{3,5\}\})$. Let $\mathcal{S} = \{P_3, K_3, K_{1,3}, C_6, C_6', C_6''\}$.

If $v \in V(G)$ such that $d_G(v) = 1$, then $v$ is called a \emph{leaf of $G$}. 
Let $L(G)$ denote the set of leaves of $G$. Let $\ell(G)$ denote $|L(G)|$. 

If $n \leq 3$, then let $B_{n,P_3}' = B_{n,P_3}$. If $n \geq 4$, then let $B_{n,P_3}'$ be the graph with $V(B_{n,P_3}') = V(B_{n,P_3})$ and $E(B_{n,P_3}') = E(B_{n,P_3}) \backslash \bigcup_{i=1}^{a_{n,3}} \{\{i,v\} \colon v \in L(F_i)\}$. Thus, if $n \geq 4$, then $B_{n,P_3}'$ is the connected $n$-vertex graph obtained from the union of $P_{a_{n,3}}, F_1, \dots, F_{a_{n,3}}$ by joining $i$ (a vertex of $P_{a_{n,3}}$) to the vertex of $F_i$ of degree $2$ for each $i \in [a_{n,3}]$, and joining $a_{n,3}$ to each $j \in [b_{n,3}] \backslash [a_{n,3}]$.

In Section~\ref{Proofsection k=2}, we prove the following result.

\begin{theorem} \label{result k=2}
If $G$ is a connected $n$-vertex graph, $G$ is not an $\mathcal{S}$-graph, and $r = \ell(G)$, then
\[\iota_2(G) \leq \left \lfloor\frac{4n-r}{14} \right \rfloor.\] 
Moreover, equality holds if $G = B_{n,P_3}'$.
\end{theorem}
It is easy to check that the inequality in Theorem~\ref{result k=2} does not hold if $G$ is an $\mathcal{S}$-graph. Note that deleting the closed neighbourhood of an $\mathcal{E}_2$-isolating set of $G$ produces a graph such that no two edges have a common vertex (and hence each component is a copy of $K_1$ or of $K_2$).

In Section~\ref{Proofsection k=3}, we prove the next result.

\begin{theorem} \label{result k=3}
If $G$ is a connected $n$-vertex graph that is neither a triangle nor a $7$-cycle, then 
\[\iota_3(G) \leq \left \lfloor\frac{n}{4} \right \rfloor.\] 
Moreover, equality holds if $G = B_{n,K_3}$.
\end{theorem}
Note that Theorem~\ref{Borgcycle} follows from Theorem~\ref{result k=3}.

We prove Theorems~\ref{result k=2} and \ref{result k=3} using the same new strategy. We build on key ideas from \cite{Borg1}, but, out of necessity, we provide a more efficient inductive argument and new ingredients, which are more abundant in the proof of Theorem~\ref{result k=2}.

\section{Proof of Theorem~\ref{result k=2}} \label{Proofsection k=2}

In this section, we prove Theorem~\ref{result k=2}. 

We start with two lemmas from \cite{Borg1} that will be used repeatedly. We provide their proofs for completeness.

\begin{lemma} \label{lemma}
If $G$ is a graph, $\mathcal{F}$ is a set of graphs, $X \subseteq V(G)$, and $Y \subseteq N[X]$, then \[\iota(G, \mathcal{F}) \leq |X| + \iota(G-Y, \mathcal{F}).\] 
\end{lemma}
\textbf{Proof.} Let $D$ be an $\mathcal{F}$-isolating set of $G-Y$ of size $\iota(G-Y, \mathcal{F})$. Clearly, $\emptyset \neq V(F) \cap Y \subseteq V(F) \cap N[X]$ for each $\mathcal{F}$-graph $F$ that is a subgraph of $G$ and not a subgraph of $G-Y$. Thus, $X \cup D$ is an $\mathcal{F}$-isolating set of $G$. The result follows.~\hfill{$\Box$}
\\

For a graph $G$, let ${\rm C}(G)$ denote the set of components of $G$.

\begin{lemma} \label{lemmacomp}
If $G$ is a graph and $\mathcal{F}$ is a set of graphs, then \[\iota(G, \mathcal{F}) = \sum_{H \in {\rm C}(G)} \iota(H, \mathcal{F}).\] 
\end{lemma}
\textbf{Proof.} For each $H \in {\rm C}(G)$, let $D_H$ be a smallest $\mathcal{F}$-isolating set of $H$. Then, $\bigcup_{H \in {\rm C}(G)} D_H$ is an $\mathcal{F}$-isolating set of $G$, so $\iota(G,\mathcal{F}) \leq \sum_{H \in {\rm C}(G)} |D_H| = \sum_{H \in {\rm C}(G)} \iota(H,\mathcal{F})$. Let $D$ be a smallest $\mathcal{F}$-isolating set of $G$. For each $H \in {\rm C}(G)$, $D \cap V(H)$ is an $\mathcal{F}$-isolating set of $H$. We have $\sum_{H \in {\rm C}(G)} \iota(H,\mathcal{F}) \leq \sum_{H \in {\rm C}(G)} |D \cap V(H)| = |D| = \iota(G,\mathcal{F})$. The result follows.~\hfill{$\Box$} \\

For a graph $G$ and a subgraph $H$ of $G$, let $L_G(H) = L(G) \cap V(H)$, $\ell_G(H) = |L_G(H)|$, and 
\[\beta_G(H) = \frac{4|V(H)| - \ell_G(H)}{14}.\] 
We may abbreviate $\beta_G(G)$ to $\beta(G)$. 

\begin{lemma} \label{betacomp} (a) If $H_1, \dots, H_t$ are pairwise vertex-disjoint subgraphs of $G$ such that $V(G) = \bigcup_{i=1}^t V(H_i)$, then 
\begin{equation} \beta(G) = \sum_{i=1}^t \beta_G(H_i). \nonumber
\end{equation}
(b) If $H$ is a connected subgraph of $G$ with $|V(H)| \geq 2$, then $\beta(H) \leq \beta_G(H)$.
\end{lemma}
\textbf{Proof.} (a) Trivially, $L(G) = \bigcup_{i=1}^t L_G(H_i)$ and
\[\beta(G) = \frac{4|V(G)|}{14} - \frac{|L(G)|}{14} = \frac{4}{14}\left( \sum_{i=1}^t |V(H_i)| \right) - \frac{1}{14}\sum_{i=1}^t |L_G(H_i)| = \sum_{i=1}^t \beta_G(H_i).\]
(b) Suppose $v \in L(G) \cap V(H)$. Since $|V(H)| \geq 2$ and $H$ is connected, $vw \in E(H)$ for some $w \in V(H)$. Since $1 \leq d_H(v) \leq d_G(v) = 1$, $d_H(v) = 1$. Thus, $L_G(H) \subseteq L(H)$, and hence $\beta(H) \leq \beta_G(H)$.~\hfill{$\Box$} \\

If $n \geq 4$, then $\{4i \colon 1 \leq i \leq n/4\}$ is an $\mathcal{E}_2$-isolating set of $P_n$, and $\{5i-4 \colon 1 \leq i \leq \lfloor (n+4)/5 \rfloor\}$ is an $\mathcal{E}_2$-isolating set of $C_n$. Consequently, 
\begin{equation} \mbox{if $n \geq 4$, $G \in \{P_n, C_n\}$, and $G \neq C_6$, then $\iota_2(G) \leq \beta(G)$.} \label{pcbound}
\end{equation}

The next lemma settles Theorem~\ref{result k=2} for $n \leq 7$.

\begin{lemma} \label{small n} If $G$ is a connected $n$-vertex graph, $n \leq 7$, and $G$ is not an $\mathcal{S}$-graph, then $\iota_2(G) \leq \lfloor \beta(G) \rfloor$.
\end{lemma}
\textbf{Proof.} Since $\iota_2(G)$ is an integer, it suffices to prove that $\iota_2(G) \leq \beta(G)$. If $n \leq 2$, then $\iota_2(G) = 0$. Since $G$ is not an $\mathcal{S}$-graph, $n \neq 3$. Suppose $4 \leq n \leq 7$. Let $v$ be a vertex of $G$ of largest degree. Since $G$ is connected, $d(v) \geq 2$ and $d(w) \geq 1$ for each $w \in V(G)$. If $d(v) = 2$, then $G$ is a path or a cycle, so $\iota_2(G) \leq \beta(G)$ by (\ref{pcbound}). Suppose $d(v) \geq 3$. Let $G' = G - N[v]$. Let $n' = |V(G')|$. Then, $n' \leq n-4$.\medskip

Suppose $n = 4$. Then, each vertex in $V(G) \backslash \{v\}$ is a neighbour of $v$. Since $G \not \simeq K_{1,3}$, $uw \in E(G)$ for some $u, w \in V(G) \backslash \{v\}$. Thus, $\ell(G) \leq 1$. Since $\{v\}$ is an $\mathcal{E}_2$-isolating set of $G$, $\iota_2(G) = 1 < \beta(G)$.\medskip  

Suppose $n \geq 5$. If no component of $G'$ has more than $2$ vertices, then $\{v\}$ is an $\mathcal{E}_2$-isolating set of $G$, so $\iota_2(G) = 1 < 15/14 \leq 3n/14 < \beta(G)$. Suppose that a component $H$ of $G'$ has at least $3$ vertices. Since $n \leq 7$ and $n' \geq |V(H)|$, we have $n = 7$, $d(v) = 3$, $n' = 3$, $G' = H$, and either $H \simeq P_3$ or $H \simeq K_3$. Since $G$ is connected, $xy \in E(H)$ for some $x \in N(v)$ and some $y \in V(H)$. Since $d(v) = 3$, $N(v) = \{x, x', x''\}$ for some $x', x'' \in N(v) \backslash \{x\}$ such that $x' \neq x''$. Since $H$ is connected and $|V(H)| = 3$, $N_H[y^*] = V(H)$ for some $y^* \in V(H)$. If $L(G) = \emptyset$, then, since $\{v, y^*\}$ is an $\mathcal{E}_2$-isolating set of $G$, $\iota_2(G) \leq 2 = \beta(G)$. Suppose $L(G) \neq \emptyset$. If we show that $G$ has an $\mathcal{E}_2$-isolating set $D$ of size $1$, then $\iota_2(G) = 1 < 3n/14 < \beta(G)$. Let $v^* \in L(G)$. Then, $v^* \notin \{v, x, y^*\}$.\medskip

Suppose $v^* \in N(v)$. Then, $N(v^*) = \{v\}$. We may assume that $v^* = x''$. If $|N(w) \cap V(H)| \geq 2$ for some $w \in \{x, x'\}$, then we take $D = \{w\}$. Suppose $|N(w) \cap V(H)| \leq 1$ for each $w \in \{x, x'\}$. Then, $N(x) \cap V(H) = \{y\}$. If $N(x') \cap V(H) = \emptyset$ or $xx' \in E(G)$, then we take $D = \{x\}$. Suppose $N(x') \cap V(H) \neq \emptyset$ and $xx' \notin E(G)$. Let $y'$ be the vertex in $N(x') \cap V(H)$. If $y' = y$, then we take $D = \{y\}$. Suppose $y' \neq y$. If we assume that $H \simeq K_3$, then we obtain $G \simeq C_6''$, a contradiction. Thus, $H \simeq P_3$. If we assume that $y^* \notin \{y, y'\}$, then we obtain $G \simeq C_6'$, a contradiction. Thus, $y^* \in \{y, y'\}$. Consequently, $y^*$ is the unique neighbour of the vertex in $V(H) \backslash \{y, y'\}$. If $y^* = y$, then we take $D = \{x\}$. If $y^* = y'$, then we take $D = \{x'\}$.\medskip

Now suppose $v^* \in V(H)$. Then, $H \simeq P_3$ and $N(v^*) = \{y^*\}$. Let $y'$ be the vertex in $V(H) \backslash \{v^*, y^*\}$. If $|N(w) \cap N(v)| \geq 2$ for some $w \in \{y^*, y'\}$, then we take $D = \{w\}$. Suppose $|N(w) \cap N(v)| \leq 1$ for each $w \in \{y^*, y'\}$. Since $x \in N(y)$, we have $y \in \{y^*, y'\}$ and $N(y) \cap N(v) = \{x\}$. If $w$ is the vertex in $\{y^*, y'\} \backslash \{y\}$ and $N(w) \cap N(v)$ is $\emptyset$ or $\{x\}$, then we take $D = \{x\}$. Now suppose $|N(y^*) \cap N(v)| = |N(y') \cap N(v)| = 1$ and $N(y^*) \cap N(v) \neq N(y') \cap N(v)$. Let $x_1 \in N(v)$ such that $\{x_1\} = N(y^*) \cap N(v)$, let $x_2 \in N(v)$ such that $\{x_2\} = N(y') \cap N(v)$, and let $x_3$ be the vertex in $N(v) \backslash \{x_1, x_2\}$. If we assume that $x_2x_3 \in E(G)$ and $x_1x_2, x_1x_3 \notin E(G)$, then we obtain $G \simeq C_6''$, a contradiction. Thus, $x_2x_3 \notin E(G)$ or $E(G) \cap \{x_1x_2, x_1x_3\} \neq \emptyset$. We take $D = \{x_1\}$.~\hfill{$\Box$}

\begin{lemma} \label{small n 2} If $G$ is an $\mathcal{S}$-graph, then
\begin{equation} \iota_2(G) \leq \beta(G) + \frac{4}{14}. \nonumber
\end{equation}
\end{lemma}
\textbf{Proof.} Trivially, $\iota_2(P_3) = 1 = \beta(P_3) + \frac{4}{14}$, $\iota_2(K_3) = 1 = \beta(K_3) + \frac{2}{14}$, $\iota_2(K_{1,3}) = 1 = \beta(K_{1,3}) + \frac{1}{14}$, $\iota_2(C_6) = 2 = \beta(C_6) + \frac{4}{14}$, $\iota_2(C_6') = 2 = \beta(C_6') + \frac{1}{14}$, and $\iota_2(C_6'') = 2 = \beta(C_6'') + \frac{1}{14}$.~\hfill{$\Box$}
\\

\noindent
\textbf{Proof of Theorem~\ref{result k=2}.} Let $B' = B_{n,P_3}'$. If $n \leq 2$, then $\iota_2(B') = \left \lfloor \beta(B') \right \rfloor = 0$. If $3 \leq n \leq 4$, then $B'$ is an $\mathcal{S}$-graph. Suppose $n \geq 5$. Since $[a_{n,3}]$ is a dominating set of $B'$, $\iota_2(B') \leq a_{n,3}$. If $D$ is an $\mathcal{E}_2$-isolating set of $B'$ of size $\iota_2(B')$, then for each $i \in [a_{n,3}]$, $D \cap (V(F_i) \cup \{i\}) \neq \emptyset$ as $B' - N_{B'}[D]$ does not contain the copy $F_i$ of $P_3$. Thus, $\iota_2(B') = a_{n,3}$. Let $c_{n,3} = b_{n,3} - a_{n,3}$. Then, $0 \leq c_{n,3} \leq 3$. Since $L(B') = ([b_{n,3}] \backslash [a_{n,3}]) \cup \bigcup_{i=1}^{a_{n,3}} L(F_i)$, $\ell(B') = c_{n,3} + 2a_{n,3}$. Since $\iota_2(B') = a_{n,3}$, 
\[ \iota_2(B') = \frac{4(4a_{n,3}) - 2a_{n,3}}{14} = \frac{4(n - c_{n,3}) + c_{n,3} - \ell(B')}{14} = \frac{4n - \ell(B') - 3c_{n,3}}{14} = \left \lfloor \beta(B') \right \rfloor .\]
We have shown that the bound in the theorem is attained if $G = B'$.\medskip

We now prove the inequality in the theorem, using induction on $n$. Since $\iota_2(G)$ is an integer, it suffices to prove that $\iota_2(G) \leq \beta(G)$. 
If $n \leq 7$, then the result is given by Lemma~\ref{small n}. Suppose $n \geq 8$.  
Let $k = \max \{d(v) \colon v \in V(G)\}$. Since $G$ is connected, $k \geq 2$. Let $v \in V(G)$ such that $d(v) = k$. If $k = 2$, then $G$ is a path or a cycle, so $\iota_2(G) \leq \beta(G)$ by (\ref{pcbound}). Suppose $k \geq 3$. Then, $|N[v]| \geq 4$. If $V(G) = N[v]$, then $\{v\}$ is an $\mathcal{E}_2$-isolating set of $G$, so $\iota_{2}(G) = 1 < 3n/14 < \beta(G)$. Suppose $V(G) \neq N[v]$. Let $G' = G-N[v]$ and $n' = |V(G')|$. Then, $n \geq n' + 4$ and $V(G') \neq \emptyset$. Let $\mathcal{H}$ be the set of components of $G'$. Let 
\begin{equation} \mathcal{H}' = \{H \in \mathcal{H} \colon \iota_2(H) > \beta_G(H)\}. \nonumber 
\end{equation} \vskip 0.1in

We have $v \notin L(G)$. Since $G$ is connected and $V(G) \neq N[v]$, $v$ has a neighbour $\overline{v}$ that has a neighbour in $V(G) \backslash N[v]$, so $\overline{v} \notin L(G)$. We have 
\begin{equation} \beta_G(G[N[v]]) \geq \frac{4|N[v]| - (|N[v]| - 2)}{14} \geq \frac{3(4)+2}{14} = 1. \label{betaN[v]}
\end{equation}  \vskip 0.1in

Suppose $\mathcal{H}' = \emptyset$. By Lemma~\ref{lemma} (with $X = \{v\}$ and $Y = N[v]$), Lemma~\ref{lemmacomp}, and Lemma~\ref{betacomp}(a),
\begin{align}
\iota_{\rm 2}(G) &\leq 1 + \iota_{\rm 2}(G') = 1 + \sum_{H \in \mathcal{H}} \iota_{\rm 2}(H) \leq \beta_G(G[N[v]]) + \sum_{H \in \mathcal{H}} \beta_G(H) = \beta(G). \nonumber 
\end{align} \vskip 0.1in

Now suppose $\mathcal{H}' \neq \emptyset$. If $H \in \mathcal{H}$ with $|V(H)| \leq 2$, then $\iota_2(H) = 0 \leq \beta_G(H)$. If $H \in \mathcal{H}$ such that $|V(H)| \geq 3$ and $H$ is not an $\mathcal{S}$-graph, then $\iota_2(H) \leq \beta(H) \leq \beta_G(H)$ by the induction hypothesis and Lemma~\ref{betacomp}(b). Thus, each member of $\mathcal{H}'$ is an $\mathcal{S}$-graph. Let $\mathcal{H}^{(1)} = \{H \in \mathcal{H}' \colon H \simeq P_3\}$, $\mathcal{H}^{(2)} = \{H \in \mathcal{H}' \colon H \simeq K_3\}$, $\mathcal{H}^{(3)} = \{H \in \mathcal{H}' \colon H \simeq C_6\}$, $\mathcal{H}^{(4)} = \{H \in \mathcal{H}' \colon H \simeq K_{1,3}, \, \ell_G(H) = 3\}$, $\mathcal{H}^{(5)} = \{H \in \mathcal{H}' \colon H \simeq C_6', \, \ell_G(H) = 1\}$, and $\mathcal{H}^{(6)} = \{H \in \mathcal{H}' \colon H \simeq C_6'', \, \ell_G(H) = 1\}$. 
If $H \in \mathcal{H}$, $H \simeq K_{1,3}$, and $xy \in E(G)$ for some $x \in N(v)$ and some $y$ in the $3$-element set $L(H)$, then $\iota_2(H) = 1 \leq \beta_G(H)$, so $H \notin \mathcal{H}'$. Similarly, if $H \in \mathcal{H}$, $H$ is a copy of $C_6'$ or of $C_6''$, $y$ is the leaf of $H$, and $xy \in E(G)$ for some $x \in N(v)$, then $\iota_2(H) = 2 = \beta_G(H)$, so $H \notin \mathcal{H}'$.  Therefore,
\begin{equation} \mathcal{H}' = \bigcup_{i=1}^6 \mathcal{H}^{(i)}. \label{calH'}
\end{equation}
Note that 
\begin{equation} \mbox{if $H \in \mathcal{H}^{(4)} \cup \mathcal{H}^{(5)} \cup \mathcal{H}^{(6)}$ and $y \in L(H)$, then $N_G(y) = N_H(y)$.} \label{calH'i}
\end{equation}  \vskip 0.1in

For any $H \in \mathcal{H}$ and any $x \in N(v)$ such that $xy_{x,H} \in E(G)$ for some $y_{x,H} \in V(H)$, we say that $H$ is \emph{linked to $x$} and that $x$ is \emph{linked to $H$}. Since $G$ is connected, each member of $\mathcal{H}$ is linked to at least one member of $N(v)$. For each $x \in N(v)$, let $\mathcal{H}_x = \{H \in \mathcal{H} \colon H \mbox{ is linked to } x\}$, $\mathcal{H}'_x = \{H \in \mathcal{H}' \colon H \mbox{ is linked to } x\}$, and $\mathcal{H}_x^* = \{H \in \mathcal{H} \backslash \mathcal{H}' \colon H \mbox{ is linked to $x$ only}\}$. For each $H \in \mathcal{H} \backslash \mathcal{H}'$, let $D_H$ be an $\mathcal{E}_2$-isolating set of $H$ of size $\iota_2(H)$.\medskip 


For any $x \in N(v)$ and any $H \in \mathcal{H}_x'$, let $H_x' = H - \{y_{x,H}\}$. By (\ref{calH'}), $H \in \mathcal{H}^{(i)}$ for some $i \in [6]$. If $i \in \{4, 5, 6\}$, then, since $x \in N(y_{x,H})$, $y_{x,H} \notin L(H)$ by (\ref{calH'i}). Thus, if $i \in \{3, 5, 6\}$, then $H$ contains a $6$-cycle $(\{y_1, \dots, y_6\}, \{y_1y_2, \dots, y_5y_6, y_6y_1\})$ 
such that $y_{x,H} = y_1$, and we take $D_{x,H} = \{y_4\}$. If $i \in \{1, 2, 4\}$, then we take $D_{x,H} = \emptyset$. If $i = 4$, then, since $y_{x,H} \notin L(H)$, $E(H_x') = \emptyset$. Clearly, $D_{x,H}$ is an $\mathcal{E}_2$-isolating set of $H_x'$, and 
\begin{equation} |D_{x,H}| \leq \beta_G(H_x') - \frac{6}{14}. \label{DxHineq}
\end{equation}
Since $\{v, y_{x,H}\} \in N(x)$, $x \notin L(G)$. Since $x \in N(y_{x,H})$ and $N_H(y_{x,H}) \neq \emptyset$, $y_{x,H} \notin L(G)$.
\\

\noindent
\emph{Case 1: $|\mathcal{H}'_x| \geq 2$ for some $x \in N(v)$.} For each $H \in \mathcal{H}' \backslash \mathcal{H}'_x$, let $x_H \in N(v)$ such that $H$ is linked to $x_H$. Let $X = \{x_H \colon H \in \mathcal{H}' \backslash \mathcal{H}'_x\}$. Note that $x \notin X$. Let 
\[D = \{v, x\} \cup X \cup \left( \bigcup_{H \in \mathcal{H}_x'} D_{x,H} \right) \cup \left( \bigcup_{H \in \mathcal{H}' \backslash \mathcal{H}_x'} D_{x_H,H} \right) \cup \left( \bigcup_{H \in \mathcal{H} \backslash \mathcal{H}'} D_{H} \right).\]
We have $V(G) = N[v] \cup \bigcup_{H \in \mathcal{H}} V(H)$, $y_{x,H} \in N[x]$ for each $H \in \mathcal{H}'_x$, and $y_{x_H,H} \in N[x_H]$ for each $H \in \mathcal{H}' \backslash \mathcal{H}'_x$, so $D$ is an $\mathcal{E}_2$-isolating set of $G$. Let 
\[Y = \{v, x\} \cup X \cup \{y_{x,H} \colon H \in \mathcal{H}_x'\} \cup \{y_{x_H,H} \colon H \in \mathcal{H}' \backslash \mathcal{H}'_x\}.\]
Then,
\begin{equation} \beta_G(G[Y]) = \frac{8 + 4|X| + 4|\mathcal{H}_x'| + 4|\mathcal{H}' \backslash \mathcal{H}_x'|}{14} \geq \frac{16 + 8|X|}{14}
 \label{GYineq}
\end{equation}
as $|\mathcal{H}'_x| \geq 2$ and $|X| \leq |\mathcal{H}' \backslash \mathcal{H}_x'|$. We have
\begin{align} \iota_2(G) &\leq |D| = 2 + |X| + \sum_{H \in \mathcal{H}_x'} |D_{x,H}| + \sum_{H \in \mathcal{H}' \backslash \mathcal{H}_x'} |D_{x_H,H}| + \sum_{H \in \mathcal{H} \backslash \mathcal{H}'} |D_{H}| \nonumber \\
&\leq 2 + |X| + \sum_{H \in \mathcal{H}_x'} \left(\beta_G(H_x') - \frac{6}{14} \right) + \sum_{H \in \mathcal{H}' \backslash \mathcal{H}_{x}'} \left( \beta_G(H_{x_H}') - \frac{6}{14} \right) \nonumber \\
& \quad + \sum_{H \in \mathcal{H} \backslash \mathcal{H}'} \beta_G(H) \quad \quad \mbox{(by (\ref{DxHineq}))} \nonumber \\
&\leq \frac{16 + 8|X|}{14} + \sum_{H \in \mathcal{H}_x'} \beta_G(H_x') + \sum_{H \in \mathcal{H}' \backslash \mathcal{H}_{x}'} \beta_G(H_{x_H}') + \sum_{H \in \mathcal{H} \backslash \mathcal{H}'} \beta_G(H) \nonumber \\ 
&\quad \; \; \mbox{(as $|\mathcal{H}'_x| \geq 2$ and $|X| \leq |\mathcal{H}' \backslash \mathcal{H}_x'|$)} \nonumber \\
&\leq \beta_G(G[Y]) + \sum_{H \in \mathcal{H}_x'} \beta_G(H_x') + \sum_{H \in \mathcal{H}' \backslash \mathcal{H}_x'} \beta_G(H_{x_H}') +  \sum_{H \in \mathcal{H} \backslash \mathcal{H}'} \beta_G(H) \quad \quad \mbox{(by (\ref{GYineq}))}\nonumber \\
&\leq \beta(G) \quad \quad \mbox{(by Lemma~\ref{betacomp}).}
\label{Hx>1}
\end{align}  \vskip 0.2in

\noindent
\emph{Case 2:}
\begin{equation} |\mathcal{H}'_x| \leq 1 \mbox{ \emph{for each} } x \in N(v). \label{Hx<2} 
\end{equation} 
For each $H \in \mathcal{H}'$, let $x_H \in N(v)$ such that $H$ is linked to $x_H$. Let $X = \{x_H \colon H \in \mathcal{H}'\}$. Then, $\mathcal{H}' = \bigcup_{x \in X} \mathcal{H}_x'$. By (\ref{Hx<2}), no two members of $\mathcal{H}'$ are linked to the same neighbour of $v$, so $|X| = |\mathcal{H}'|$. 
Let $W = N(v) \backslash X$, $W' = N[v] \backslash X$, and $Y_X = X \cup \{y_{x_H,H} \colon H \in \mathcal{H}'\}$. We have $|Y_X| = |X| + |\mathcal{H}'| = 2|X|$ and $\beta_G(G[Y_X]) = \frac{8|X|}{14}$. Let 
\[D = \{v\} \cup X \cup \left( \bigcup_{H \in \mathcal{H}'}D_{x_H,H} \right) \cup \left( \bigcup_{H \in \mathcal{H} \backslash \mathcal{H}'} D_{H} \right).\] 
We have $V(G) = N[v] \cup \bigcup_{H \in \mathcal{H}} V(H)$ and $y_{x_H,H} \in N[x_H]$ for each $H \in \mathcal{H}'$, so $D$ is an $\mathcal{E}_2$-isolating set of $G$.\medskip

\noindent
\emph{Subcase 2.1: $|W| \geq 3$.} If $\beta_G(G[W]) \geq \frac{10}{14}$, then $\beta_G(G[W']) \geq 1$ and, similarly to (\ref{Hx>1}),
\begin{align} \iota_2(G) &\leq |D| = 1 + |X| + \sum_{H \in \mathcal{H}'} |D_{x_H,H}| + \sum_{H \in \mathcal{H} \backslash \mathcal{H}'} |D_{H}| \nonumber \\
&\leq \beta_G(G[W']) + |X| + \sum_{H \in \mathcal{H}'} \left(\beta_G(H_{x_H}') - \frac{6}{14} \right) + \sum_{H \in \mathcal{H} \backslash \mathcal{H}'} \beta_G(H) \nonumber \\
&= \beta_G(G[W']) + \frac{8|X|}{14} + \sum_{H \in \mathcal{H}'} \beta_G(H_{x_H}') + \sum_{H \in \mathcal{H} \backslash \mathcal{H}'} \beta_G(H) \nonumber \\
&= \beta_G(G[W']) + \beta_G(G[Y_X]) + \sum_{H \in \mathcal{H}'} \beta_G(H_{x_H}') + \sum_{H \in \mathcal{H} \backslash \mathcal{H}'} \beta_G(H) = \beta(G). \nonumber
\end{align}
Suppose $\beta_G(G[W]) < \frac{10}{14}$. We have $\frac{9}{14} \geq \beta_G(G[W]) = \frac{4|W| - \ell_G(G[W])}{14} \geq \frac{3|W|}{14} \geq \frac{9}{14}$, so $|W| = 3$ and $W \subseteq L(G)$. Thus, since $v \in N[X]$, we have that $D \backslash \{v\}$ is an $\mathcal{E}_2$-isolating set of $G$, and hence, similarly to the above, 
\begin{align} \iota_2(G) &\leq |D| - 1 \leq \beta_G(G[Y_X]) + \sum_{H \in \mathcal{H}'} \beta_G(H_{x_H}') + \sum_{H \in \mathcal{H} \backslash \mathcal{H}'} \beta_G(H) \nonumber \\
&< \beta_G(G[W']) + \beta_G(G[Y_X]) + \sum_{H \in \mathcal{H}'} \beta_G(H_{x_H}') + \sum_{H \in \mathcal{H} \backslash \mathcal{H}'} \beta_G(H) = \beta(G).\nonumber
\end{align} \vskip 0.1in

\noindent
\emph{Subcase 2.2: $|W| \leq 2$.}\medskip

\noindent
\emph{Subsubcase 2.2.1: Some member $H$ of $\mathcal{H}'$ is linked to $x_H$ only.} Let $G^* = G - (\{x_H\} \cup V(H))$. Then, $G^*$ has a component $G_v^*$ such that $N[v] \backslash \{x_H\} \subseteq V(G_v^*)$, and the other components of $G^*$ are the members of $\mathcal{H}_{x_H}^*$. Let $D^*$ be an $\mathcal{E}_2$-isolating set of $G_v^*$ of size $\iota_2(G_v^*)$. Let $D' = D^* \cup \{x_H\} \cup D_{x_H,H} \cup \bigcup_{I \in \mathcal{H}_{x_H}^*} D_I$. Since $D'$ is an $\mathcal{E}_2$-isolating set of $G$,
\begin{align} \iota_2(G) &\leq |D^*| + 1 + |D_{x_H,H}| + \sum_{I \in \mathcal{H}_{x_H}^*} |D_I| \nonumber \\
&\leq \iota_2(G_v^*) + \beta_G(G[\{x_H,y_{x_H,H}\}]) + \beta_G(H_x') + \sum_{I \in \mathcal{H}_{x_H}^*} \beta_G(I) \quad \mbox{(by \ref{DxHineq})}.  \nonumber
\end{align}
If $\iota_2(G_v^*) \leq \beta_G(G_v^*)$, then $\iota_2(G) \leq \beta(G)$ by Lemma~\ref{betacomp}(a). Suppose $\iota_2(G_v^*) > \beta_G(G_v^*)$. Since $|V(G_v^*)| \geq |N[v] \backslash \{x_H\}| \geq 3$, $\beta(G_v^*) \leq \beta_G(G_v^*)$ by Lemma~\ref{betacomp}(b). Thus, $\iota_2(G_v^*) > \beta(G_v^*)$. By the induction hypothesis, $G_v^*$ is an $\mathcal{S}$-graph. We have $v \in N[x_H]$ and $d_{G_v^*}(v) \geq 2$ (so $v \notin L(G_v^*)$). Thus, if $G_v^*$ is a $\{P_3, K_3, K_{1,3}\}$-graph, then $D' \backslash D^*$ is an $\mathcal{E}_2$-isolating set of $G$, and hence
\[\iota_2(G) \leq \beta_G(G[\{x_H,y_{x_H,H}\}]) + \beta_G(H_x') + \sum_{I \in \mathcal{H}_{x_H}^*} \beta_G(I) < \beta(G)\] 
by Lemma~\ref{betacomp}(a). Suppose that $G_v^*$ is a $\{C_6, C_6', C_6''\}$-graph. Then, $G_v^*$ contains a $6$-cycle $(\{y_1, \dots, y_6\}, \{y_1y_2, \dots, y_5y_6, y_6y_1\})$ such that $v = y_1$ (recall that $v \notin L(G_v^*)$). Thus, since $v \in N[x_H]$, $(D' \backslash D^*) \cup \{y_4\}$ is an $\mathcal{E}_2$-isolating set of $G$, and hence 
\begin{align} \iota_2(G) &\leq 1 + \beta_G(G[\{x_H,y_{x_H,H}\}]) + \beta_G(H_x') + \sum_{I \in \mathcal{H}_{x_H}^*} \beta_G(I) \nonumber \\
&< \beta_G(G_v^*) + \beta_G(G[\{x_H,y_{x_H,H}\}]) + \beta_G(H_x') + \sum_{I \in \mathcal{H}_{x_H}^*} \beta_G(I) = \beta(G) \nonumber
\end{align}
by Lemma~\ref{betacomp}(a).\medskip

\noindent
\emph{Subsubcase 2.2.2: For each $H \in \mathcal{H}'$, $H$ is linked to some $x_H' \in N(v) \backslash \{x_H\}$.} Let $X' = \{x_H' \colon H \in \mathcal{H}'\}$. By (\ref{Hx<2}), for every $H, I \in \mathcal{H}'$ with $H \neq I$, we have $x_H' \neq x_I'$ and $x_H \neq x_I'$. Thus, $|X'| = |\mathcal{H}'|$ and $X' \subseteq W$. This yields $|\mathcal{H}'| \leq 2$ as $|W| \leq 2$.\medskip

Suppose $|\mathcal{H}'| = 2$. Then, $|W| = 2$ and $X' = W$. Let $H_1$ and $H_2$ be the two members of $\mathcal{H}'$. Let $x_1 = x_{H_1}$, $x_2 = x_{H_2}$, $x_1' = x_{H_1}'$, and $x_2' = x_{H_2}'$. We have $N(v) = X \cup W = \{x_1, x_2, x_1', x_2'\}$. Thus, $N(v) \cap L(G) = \emptyset$. By (\ref{calH'}), $H_1, H_2 \in \bigcup_{i=1}^6 \mathcal{H}^{(i)}$.\medskip 

Suppose $H_1 \in \mathcal{H}^{(2)} \cup \mathcal{H}^{(4)} \cup \mathcal{H}^{(5)} \cup \mathcal{H}^{(6)}$. Let $Y_1 = \{x_1\} \cup V(H_1)$ and $D_1 = \{y_{x_1, H_1}\} \cup D_{x_1, H_1}$. Clearly, $|D_1| \leq \beta_G(G[Y_1]) - \frac{2}{14}$. Let $G^* = G - Y_1$. Then, $G^*$ has a component $G_v^*$ such that $N[v] \backslash \{x_1\} \subseteq V(G_v^*)$, and the other components of $G^*$ are the members of $\mathcal{H}_{x_1}^*$. Since $(N[v] \backslash \{x_1\}) \cup V(H_2) \subseteq V(G_v^*)$, $|V(G_v^*)| \geq 7$. Let $D^*$ be an $\mathcal{E}_2$-isolating set of $G_v^*$ of size $\iota_2(G_v^*)$. If $G_v^*$ is an $\mathcal{S}$-graph, then $G_v^*$ is a $\{C_6', C_6''\}$-graph, so $\iota_2(G_v^*) = 2 \leq \beta_G(G_v^*) + \frac{1}{14}$. If $G_v^*$ is not an $\mathcal{S}$-graph, then $\iota_2(G_v^*) \leq \beta(G_v^*) \leq \beta_G(G_v^*)$ by the induction hypothesis and Lemma~\ref{betacomp}(b). Note that, since $H_1 \in \mathcal{H}^{(2)} \cup \mathcal{H}^{(4)} \cup \mathcal{H}^{(5)} \cup \mathcal{H}^{(6)}$, if $y \in Y_1$ such that $y$ is a vertex of $G - N[D_1]$, then $H_1 \in \mathcal{H}^{(5)} \cup \mathcal{H}^{(6)}$, $y$ is the leaf of $H_1$, and $N_{G-N[D_1]}(y) = \emptyset$. Thus, $D^* \cup D_1 \cup \bigcup_{H \in \mathcal{H}_{x_1}^*} D_H$ is an $\mathcal{E}_2$-isolating set of $G$, and hence
\begin{align} \iota_2(G) &\leq |D^*| + |D_1| + \sum_{H \in \mathcal{H}_{x_1}^*} |D_H| \nonumber \\
&\leq \beta_G(G_v^*) + \frac{1}{14} + \beta_G(G[Y_1]) - \frac{2}{14} + \sum_{H \in \mathcal{H}_{x_1}^*} \beta_G(H) < \beta(G)   \nonumber
\end{align}
by Lemma~\ref{betacomp}(a). Similarly, $\iota_2(G) < \beta(G)$ if $H_2 \in \mathcal{H}^{(2)} \cup \mathcal{H}^{(4)} \cup \mathcal{H}^{(5)} \cup \mathcal{H}^{(6)}$.\medskip

Now suppose $H_1, H_2 \in \mathcal{H}^{(1)} \cup \mathcal{H}^{(3)}$. For any $H \simeq P_3$, let $y_H^1$ and $y_H^3$ be the two leaves of $H$, and let $y_H^2$ be the remaining vertex of $H$ (so $N_H[y_H^2] = V(H)$). If $H_1, H_2 \in \mathcal{H}^{(3)}$, then by applying the argument for the previous case, we obtain that $|D_1| = \beta_G(G[Y_1])$ ($= 2$) and that $G_v^*$ is not an $\mathcal{S}$-graph (as $|V(G_v^*)| > 7$), and $\iota_2(G) \leq \beta(G)$ follows. Suppose $H_1 \in \mathcal{H}^{(1)}$ and $H_2 \in \mathcal{H}^{(3)}$. If $y_{H_1}^1, y_{H_1}^3 \in L(G)$, then $y_{H_1}^2 = y_{x_1,H_1}$, and $\iota_2(G) \leq \beta(G)$ follows as in the case $H_1, H_2 \in \mathcal{H}^{(3)}$. If $y_{H_1}^1, y_{H_1}^3 \notin L(G)$, then $\{v, x_1, x_2\} \cup D_{H_2} \cup \bigcup_{H \in \mathcal{H} \backslash \mathcal{H}'} D_H$ is an $\mathcal{E}_2$-isolating set of $G$, so
\begin{equation} \iota_2(G) \leq 4 + \sum_{H \in \mathcal{H} \backslash \mathcal{H}'} |D_H| = \beta_G(G[N[v] \cup V(H_1) \cup V(H_2)]) + \sum_{H \in \mathcal{H} \backslash \mathcal{H}'} \beta_G(H) = \beta(G)  \nonumber
\end{equation}
by Lemma~\ref{betacomp}(a). Suppose that for some $i \in \{1, 3\}$, $y_{H_1}^i \notin L(G)$ and $y_{H_1}^{4-i} \in L(G)$. Then, $xy_{H_1}^i \in E(G)$ for some $x \in N(v)$. By (\ref{Hx<2}), $x \neq x_2$. Let $Y_1 = \{x\} \cup V(H_1)$. Let $G^* = G - Y_1$. As above, $G^*$ has a component $G_v^*$ such that $N[v] \backslash \{x\} \subseteq V(G_v^*)$, and the other components of $G^*$ are the members of $\mathcal{H}_{x}^*$. Let $D^*$ be an $\mathcal{E}_2$-isolating set of $G_v^*$ of size $\iota_2(G_v^*)$. Since $|V(G_v^*)| \geq |N[v] \backslash \{x\}| + |V(H_2)| > 7$, $G_v^*$ is not an $\mathcal{S}$-graph, so $\iota_2(G_v^*) \leq \beta(G_v^*) \leq \beta_G(G_v^*)$ by the induction hypothesis and Lemma~\ref{betacomp}(b). Since $Y_1 \backslash \{y_{H_1}^{4-i}\} \subseteq N[y_{H_1}^i]$ and $y_{H_1}^{4-i} \in L(G)$, $\{y_{H_1}^i\} \cup D^* \cup \bigcup_{H \in \mathcal{H}_{x}^*} D_H$ is an $\mathcal{E}_2$-isolating set of $G$, so
\begin{equation} \iota_2(G) \leq 1 + |D^*| + \sum_{H \in \mathcal{H}_{x}^*} |D_H| < \beta_G(G[Y_1]) + \beta_G(G_v^*) + \sum_{H \in \mathcal{H}_{x}^*} \beta_G(H) \leq \beta(G)  \nonumber
\end{equation}
by Lemma~\ref{betacomp}(a). Similarly, $\iota_2(G) < \beta(G)$ if $H_1 \in \mathcal{H}^{(3)}$ and $H_2 \in \mathcal{H}^{(1)}$. Now suppose $H_1, H_2 \in \mathcal{H}^{(1)}$. Let $Y = N[v] \cup V(H_1) \cup V(H_2)$. We have 
\[V(G) = Y \cup \bigcup_{H \in \mathcal{H} \backslash \mathcal{H}'} V(H) = \{v, x_1, x_2, x_1', x_2'\} \cup \bigcup_{i=1}^2 \{y_{H_i}^1, y_{H_i}^2, y_{H_i}^3\} \cup \bigcup_{H \in \mathcal{H} \backslash \mathcal{H}'} V(H),\]
$N[v] \cap L(G) = \emptyset$, and $d_{H_1}(y_{H_1}^2) = d_{H_2}(y_{H_2}^2) = 2$ (so $y_{H_1}^2, y_{H_2}^2 \notin L(G)$). Let $h = |\{y_{H_1}^1, y_{H_1}^3, y_{H_2}^1, y_{H_2}^3\} \cap L(G)|$. We have $\beta_G(G[Y]) = \frac{4(11) - h}{14}$. Since $\{v, y_{H_1}^2, y_{H_2}^2\}$ is a dominating set of $G[Y]$, $\{v, y_{H_1}^2, y_{H_2}^2\} \cup \bigcup_{H \in \mathcal{H} \backslash \mathcal{H}'} D_H$ is an $\mathcal{E}_2$-isolating set of $G$, so $\iota_2(G) \leq 3 + \sum_{H \in \mathcal{H} \backslash \mathcal{H}'} \beta_G(H)$. If $h \leq 2$, then $\iota_2(G) \leq \beta_G(G[Y]) + \sum_{H \in \mathcal{H} \backslash \mathcal{H}'} \beta_G(H) = \beta(G)$ by Lemma~\ref{betacomp}(a). Suppose $h \geq 3$. Then, $|L_G(H_i)| = 2$ for some $i \in [2]$, and $|L_G(H_{3-i})| \geq 1$. We may assume that $i = 2$. Thus, $y_{H_2}^1, y_{H_2}^3 \in L(G)$. Also, $y_{H_1}^j \in L(G)$ for some $j \in \{1, 3\}$, so $y_{x_1,H_1} \neq y_{H_1}^j$. Let $Y_1 = \{x_1\} \cup (V(H_1) \backslash \{y_{H_1}^j\})$ and $D_1 = \{y_{x_1,H_1}\}$. We have $Y_1 \subseteq N[D_1]$ and $N_{G-Y_1}(y_{H_1}^j) = \emptyset$. Let $G^* = G - (Y_1 \cup \{y_{H_1}^j\})$. Then, $G^*$ has a component $G_v^*$ such that $N[v] \backslash \{x_1\} \subseteq V(G_v^*)$, and the other components of $G^*$ are the members of $\mathcal{H}_{x_1}^*$. Since $(N[v] \backslash \{x_1\}) \cup V(H_2) \subseteq V(G_v^*)$, $|V(G_v^*)| \geq 7$. Since $y_{H_2}^1, y_{H_2}^3 \in L(G_v^*)$, $G_v^*$ is not an $\mathcal{S}$-graph, so $\iota_2(G_v^*) \leq \beta(G_v^*) \leq \beta_G(G_v^*)$ by the induction hypothesis and Lemma~\ref{betacomp}(b). Let $D^*$ be an $\mathcal{E}_2$-isolating set of $G_v^*$ of size $\iota_2(G_v^*)$. Then, $D^* \cup D_1 \cup \bigcup_{H \in \mathcal{H}_{x_1}^*} D_H$ is an $\mathcal{E}_2$-isolating set of $G$, and hence
\begin{align} \iota_2(G) &\leq |D^*| + |D_1| + \sum_{H \in \mathcal{H}_{x_1}^*} |D_H| \nonumber \\
&\leq \beta_G(G_v^*) + \beta_G(G[Y_1 \cup \{y_{H_1}^j\}]) + \sum_{H \in \mathcal{H}_{x_1}^*} \beta_G(H) = \beta(G)  \quad \mbox{(by Lemma~\ref{betacomp}(a)).} \nonumber
\end{align}
\vskip 0.1in   

Finally, suppose $|\mathcal{H}'| = 1$. Let $H_1$ be the member of $\mathcal{H}'$, and let $x_1 = x_{H_1}$ and $x_1' = x_{H_1}'$. We have $3 \leq d(v) = |X| + |W| = 1 + |W| \leq 3$, so $d(v) = 3$, and hence $N(v) = \{x_1, x_1', w\}$ for some $w \in V(G) \backslash \{x_1, x_1'\}$. Let $Y = N[v] \cup V(H_1)$. Then, 
\begin{equation} V(G) = V(G[Y]) \cup \bigcup_{H \in \mathcal{H} \backslash \mathcal{H}'} V(H). \label{GY}
\end{equation}
By (\ref{calH'}), $H_1 \in \mathcal{H}^{(j)}$ for some $j \in [6]$.\medskip

Suppose $j \in \{4, 5, 6\}$. Let $D_1 = \{v, y_{x_1,H_1}\} \cup D_{x_1,H_1}$. Then, $D_1 \cup \bigcup_{H \in \mathcal{H} \backslash \mathcal{H}'} D_H$ is an $\mathcal{E}_2$-isolating set of $G$, and since $v, x_1, x_1' \notin L(G)$, we clearly have $|D_1| \leq \beta_G(G[Y])$. Thus, $\iota_2(G) \leq \beta_G(G[Y]) + \sum_{H \in \mathcal{H} \backslash \mathcal{H}'} \beta_G(H) = \beta(G)$ by Lemma~\ref{betacomp}(a).\medskip

Suppose $j = 3$. We first show that $G[Y]$ has an $\mathcal{E}_2$-isolating set $D_Y$ of size~$2$. We have $H_1 \simeq C_6$. We may assume that $H_1 = C_6$. For each $i \in [6]$, let $y_i$ be the vertex $i$ of $C_6$. We may assume that $y_{x_1,H_1} = y_1$. Let $Y_1 = \{x_1, y_1, y_2, y_6\}$. We have $Y_1 \subseteq N[y_1]$. Let $I = G[Y] - Y_1$. Then, $V(I) = \{v, x_1', w, y_3, y_4, y_5\}$ and $vx_1', vw, y_3y_4, y_4y_5 \in E(I)$. If $I$ is connected and $I \not\simeq C_6$, then, by Lemma~\ref{small n}, $I$ has an $\mathcal{E}_2$-isolating set $D_{I}$ of size $1$, so we take $D_Y = \{y_1\} \cup D_{I}$. Suppose that $I$ is not connected. Then, $E(I) \backslash \{x_1'w\} = \{vx_1', vw, y_3y_4, y_4y_5\}$, and hence $N(y_i) \cap \{x_1', w\} = \emptyset$ for each $i \in \{3, 4, 5\}$. If $N(y) \cap \{x_1', w\} \neq \emptyset$ for some $y \in \{y_2, y_6\}$, then we take $D_Y = \{y_1, y\}$. If $N(y_2) \cap \{x_1', w\} = N(y_6) \cap \{x_1', w\} = \emptyset$, then we take $D_Y = \{x_1, y_4\}$. Now suppose $I \simeq C_6$. We may assume that $E(I) = \{vx_1', x_1'y_3, y_3y_4, y_4y_5, y_5w, wv\}$. Let $I' = Y - N_{G[Y]}[\{x_1, y_4\}]$. Then, $V(I') \subseteq \{x_1', w, y_2, y_6\}$. Since $G[Y] - N[v] = H_1 = C_6$, $y_2y_6 \notin E(G)$. If $uy \notin E(G)$ for each $u \in \{x_1', w\}$ and each $y \in \{y_2, y_6\}$, then $I'$ contains no $3$-path, so we can take $D_Y = \{x_1, y_4\}$. Suppose $uy \in E(G)$ for some $u \in \{x_1', w\}$ and $y \in \{y_2, y_6\}$. Since $k = d(v) = 3$, we have $d(u) = d(y) = 3$, $x_1'w \notin E(G)$, and $uy', yu' \notin E(G)$ for $y' \in \{y_2, y_6\} \backslash \{y\}$ and $u' \in \{x_1', w\} \backslash \{u\}$. Thus, $I'$ contains no $3$-path, and hence again we can take $D_Y = \{x_1, y_4\}$.\medskip

Therefore, $\iota_2(G[Y]) = 2$. Since $L_G(G[Y]) \subseteq \{w\}$, $\beta_G(G[Y]) \geq \frac{4(10) - 1}{14} = \frac{39}{14}$. If $V(G) = Y$, then $G = G[Y]$, so $\iota_2(G) = 2 < \beta(G)$. Suppose $V(G) \neq Y$. Then, $\mathcal{H} \backslash \mathcal{H}' \neq \emptyset$. Thus, $uz \in E(G)$ for some $u \in N(v)$ and some $z \in V(G) \backslash Y$. 
Let $u' \in \{x_1, x_1'\} \backslash \{u\}$. Suppose $u' = x_1$.  Suppose that $I$ is connected. Let $G^* = G - Y_1$. Then, $G^*$ has a component $G_v^*$ such that $V(I) \cup \{z\} \subseteq V(G_v^*)$, and the other components of $G^*$ are the members of $\mathcal{H}_{x_1}^*$. We have $|V(G_v^*)| \geq 7$. If $G_v^*$ is a $\{C_6', C_6''\}$-graph, then $\iota_2(G_v^*) = 2 \leq \beta_G(G_v^*) + \frac{1}{14}$; otherwise, $\iota_2(G_v^*) \leq \beta(G_v^*) \leq \beta_G(G_v^*)$ by the induction hypothesis and Lemma~\ref{betacomp}(b). Let $D^*$ be an $\mathcal{E}_2$-isolating set of $G_v^*$ of size $\iota_2(G_v^*)$. Then, $D^* \cup \{y_1\} \cup \bigcup_{H \in \mathcal{H}_{x_1}^*} D_H$ is an $\mathcal{E}_2$-isolating set of $G$, and hence
\begin{align} \iota_2(G) &\leq |D^*| + 1 + \sum_{H \in \mathcal{H}_{x_1}^*} |D_H| \nonumber \\
&\leq \left( \beta_G(G_v^*)  + \frac{1}{14} \right) + \left( \beta_G(G[Y_1]) - \frac{2}{14} \right) + \sum_{H \in \mathcal{H}_{x_1}^*} \beta_G(H) = \beta(G) - \frac{1}{14}  \nonumber
\end{align}
by Lemma~\ref{betacomp}(a). Now suppose that $I$ is not connected. Then, $xy_i \notin E(G)$ for each $x \in \{x_1', w\}$ and each $i \in \{3, 4, 5\}$. Thus, $y_{x_1',H_1} \in \{y_1, y_2, y_6\}$. Since $Y_1 \subseteq N[y_1]$, $d(y_1) \geq 3$. Since $k = d(v) = 3$, we obtain $d(y_1) = 3$, so $N[y_1] = Y_1$. Thus, $x_1'y_1 \notin E(G)$, and hence $y_{x_1',H_1} \in \{y_2, y_6\}$. Let $G^* = G - (Y \backslash \{v, w\})$. Then, $G^*$ has a component $G_v^*$ such that $\{v, w\} \subseteq V(G_v^*)$, and the other components of $G^*$ are the members of $\mathcal{H} \backslash (\mathcal{H}' \cup \mathcal{H}_w)$. If $G_v^*$ is an $\mathcal{S}$-graph, then $\iota_2(G_v^*) \leq \beta(G_v^*) + \frac{4}{14} \leq \beta_G(G_v^*) + \frac{4}{14}$ by Lemmas~\ref{small n 2} and \ref{betacomp}(b); otherwise, $\iota_2(G_v^*) \leq \beta(G_v^*) \leq \beta_G(G_v^*)$ by the induction hypothesis and Lemma~\ref{betacomp}(b). Let $D^*$ be an $\mathcal{E}_2$-isolating set of $G_v^*$ of size $\iota_2(G_v^*)$. Then, $D^* \cup \{y_1, y_{x_1',H_1}\} \cup \bigcup_{H \in \mathcal{H} \backslash (\mathcal{H}' \cup \mathcal{H}_w)} D_H$ is an $\mathcal{E}_2$-isolating set of $G$, and hence
\begin{align} \iota_2(G) &\leq |D^*| + 2 + \sum_{H \in \mathcal{H} \backslash (\mathcal{H}' \cup \mathcal{H}_w)} |D_H| \nonumber \\
&\leq \left( \beta_G(G_v^*)  + \frac{4}{14} \right) + \left( \beta_G(G[Y \backslash \{v, w\}]) - \frac{4}{14} \right) + \sum_{H \in \mathcal{H} \backslash (\mathcal{H}' \cup \mathcal{H}_w)} \beta_G(H) = \beta(G) \nonumber
\end{align}
by Lemma~\ref{betacomp}(a). Similarly, $\iota_2(G) \leq \beta(G)$ if $u' = x_1'$.\medskip

Now suppose $j \in \{1, 2\}$. Then, $H_1 \simeq P_3$ or $H_1 \simeq K_3$. We conclude by applying arguments similar to those used above. We have $Y = \{v, x_1, x_1', w, y_1, y_2, y_3\}$, where $\{y_1, y_2, y_3\} = V(H_1)$ and $N_{H_1}[y_2] = V(H_1)$. Also, $L_G(G[Y]) \subseteq \{w, y_1, y_3\}$. Since $n \geq 8$, $uz \in E(G)$ for some $u \in N(v)$ and some $z \in V(G) \backslash Y$. Let $u' \in \{x_1, x_1'\} \backslash \{u\}$. As above, we may assume that $u' = x_1$. If $L_G(G[Y]) = \emptyset$, then, since $\{v, y_2\} \cup \bigcup_{H \in \mathcal{H} \backslash \mathcal{H}'} D_H$ is an $\mathcal{E}_2$-isolating set of $G$, we have $\iota_2(G) \leq 2 + \sum_{H \in \mathcal{H} \backslash \mathcal{H}'} |D_H| = \beta_G(G[Y]) + \sum_{H \in \mathcal{H} \backslash \mathcal{H}'} \beta_G(H) = \beta(G)$ by Lemma~\ref{betacomp}(a). Suppose $L_G(G[Y]) \neq \emptyset$.\medskip

Suppose $w \in L(G)$. Then, $u = x_1'$. Let $Y_1 = \{v, x_1, y_{x_1,H_1}, w\}$. Then, $Y_1 \backslash \{w\} \subseteq N[x_1]$ and $d_{G - (Y_1 \backslash \{w\})}(w) = 0$. Let $y_1'$ and $y_2'$ be the two vertices in $V(H_1) \backslash \{y_{x_1,H_1}\}$. Let $G^* = G-Y_1$. Then, $G^*$ has a component $G_v^*$ such that $\{x_1', z\} \subseteq V(G_v^*)$, the members of $\mathcal{H}_{x_1}^*$ are components of $G^*$, and the vertex set of any other component of $G^*$ is a subset of $\{y_1', y_2'\} \backslash V(G_v^*)$. Let $D^*$ be an $\mathcal{E}_2$-isolating set of $G_v^*$ of size $\iota_2(G_v^*)$. Then, $D^* \cup \{x_1\} \cup \bigcup_{H \in \mathcal{H}_{x_1}^*} D_H$ is an $\mathcal{E}_2$-isolating set of $G$. Let $I = G[V(G_v^*) \cup \{y_1', y_2'\}]$. If $y_1', y_2' \in V(G_v^*)$, then $G_v^*$ is not a $(\mathcal{S} \backslash \{K_{1,3}\})$-graph, so we have $|D^*| \leq \beta_G(G_v^*) + \frac{1}{14} = \beta_G(I) + \frac{1}{14}$ by the induction hypothesis and the equality $\iota_2(K_{1,3}) = \beta(K_{1,3}) + \frac{1}{14}$. If $y_i' \notin V(G_v^*)$ for some $i \in \{1, 2\}$, then, by the induction hypothesis and Lemma~\ref{small n 2}, $|D^*| \leq \beta_G(G_v^*) + \frac{4}{14} \leq \beta_G(G_v^*) + \beta_G(G[\{y_i'\}]) + \frac{1}{14} \leq \beta_G(I) + \frac{1}{14}$. Since $Y_1 \cap L(G) = \{w\}$, $\beta_G(G[Y_1]) = \frac{4(4)-1}{14} = 1 + \frac{1}{14}$. We have
\begin{align} \iota_2(G) &\leq |D^*| + 1 + \sum_{H \in \mathcal{H}_{x_1}^*} |D_H| \nonumber \\
&\leq \left( \beta_G(I)  + \frac{1}{14} \right) + \left( \beta_G(G[Y_1]) - \frac{1}{14} \right) + \sum_{H \in \mathcal{H}_{x_1}^*} \beta_G(H) = \beta(G). \nonumber
\end{align}
\vskip 0.1in

Now suppose $w \notin L(G)$. Since $\emptyset \neq L_G(G[Y]) \subseteq \{w, y_1, y_3\}$, we have $y_1 \in L(G)$ or $y_3 \in L(G)$. We may assume that $y_3 \in L(G)$. Let $y^* = y_{x_1,H_1}$. Suppose $y^* \neq y_{x_1',H_1}$. Then, $\{y^*, y_{x_1',H_1}\} = \{y_1, y_2\}$, and hence $y_1, y_2 \notin L(G)$. Let $Y_1 = \{x_1\} \cup V(H_1)$. Then, $Y_1 \backslash \{y_3\} \subseteq N[y^*]$ and $d_{G - (Y_1 \backslash \{y_3\})}(y_3) = 0$. Let $G^* = G-Y_1$. Then, $G^*$ has a component $G_v^*$ such that $\{v, x_1', w, z\} \subseteq V(G_v^*)$, and the other components of $G^*$ are the members of $\mathcal{H}_{x_1}^*$. Let $D^*$ be an $\mathcal{E}_2$-isolating set of $G_v^*$ of size $\iota_2(G_v^*)$. Then, $D^* \cup \{y^*\} \cup \bigcup_{H \in \mathcal{H}_{x_1}^*} D_H$ is an $\mathcal{E}_2$-isolating set of $G$. Since $|V(G_v^*)| \geq 4$, $G_v^*$ is not a $\{P_3, K_3\}$-graph. If we assume that $G_v^* \simeq C_6$, then we obtain $G[V(G_v^*) \backslash N[v]] \simeq P_3$, which contradicts $\mathcal{H}' = \{H_1\}$. Therefore, if $G_v^*$ is an $\mathcal{S}$-graph, then $G_v^*$ is a $\{K_{1,3}, C_6', C_6''\}$-graph, and hence $|D^*| \leq \beta_G(G_v^*) + \frac{1}{4}$. Since $Y_1 \cap L(G) = \{y_3\}$, $\beta_G(G[Y_1]) = 1 + \frac{1}{14}$. As above, this yields $\iota_2(G) \leq \beta(G)$. Now suppose $y^* = y_{x_1',H_1}$. Then, $x_1, x_1' \in N(y^*)$. Let $Y_1 = \{x_1, x_1'\} \cup V(H_1)$. Then, $Y_1 \backslash \{y_3\} \subseteq N[y^*]$ and $d_{G - (Y_1 \backslash \{y_3\})}(y_3) = 0$. Let $G^* = G-Y_1$. Then, $G^*$ has a component $G_v^*$ such that $\{v, w\} \subseteq V(G_v^*)$, and the other components of $G^*$ are the members of $\mathcal{H} \backslash (\mathcal{H}' \cup \mathcal{H}_w)$. As above, $\iota_2(G_v^*) \leq \beta_G(G_v^*) + \frac{4}{14}$. Since $x_1, x_1', y^* \notin L(G)$, $\beta_G(G[Y_1]) \geq \frac{4(5) - 2}{14} = 1 + \frac{4}{14}$. Let $D^*$ be an $\mathcal{E}_2$-isolating set of $G_v^*$ of size $\iota_2(G_v^*)$. Then, $D^* \cup \{y^*\} \cup \bigcup_{H \in \mathcal{H} \backslash (\mathcal{H}' \cup \mathcal{H}_w)} D_H$ is an $\mathcal{E}_2$-isolating set of $G$, and hence $\iota_2(G) \leq \beta(G)$ as above.~\hfill{$\Box$}

\section{Proof of Theorem~\ref{result k=3}} \label{Proofsection k=3}

In this section, we prove Theorem~\ref{result k=3}, using the line of argument in the proof of Theorem~\ref{result k=2}. Again, we first settle the case $n \leq 7$.

\begin{lemma} \label{small n k=3} If $G$ is a connected $n$-vertex graph, $n \leq 7$, and $G$ is not a $\{K_3, C_7\}$-graph, then $\iota_3(G) \leq \lfloor n/4 \rfloor$.
\end{lemma}
\textbf{Proof.} If $n \leq 3$, then $|E(G)| \leq 2$ unless $G \simeq K_3$. Suppose $4 \leq n \leq 7$. Let $k = \max\{d(v) \colon v \in V(G)\}$. Since $G$ is connected, $k \geq 2$. If $k = 2$, then $G$ is an $n$-path or an $n$-cycle, and hence we clearly have $\iota_3(G) = 1 \leq \lfloor n/4 \rfloor$ unless $G \simeq C_7$. Suppose $k \geq 3$. Let $v \in V(G)$ with $d(v) = k$, $G' = G - N[v]$, $X = N(v)$, and $Y = V(G')$. Then, $|E(G')| \leq 2$ unless $k = 3$, $n = 7$, and $G' \simeq K_3$. Suppose $k = 3$, $n = 7$, and $G' \simeq K_3$. If we show that $G$ has an $\mathcal{E}_3$-isolating set $D$ of size $1$, then $\iota_3(G) \leq \lfloor n/4 \rfloor$. For each $y \in Y$, $|N(y) \cap X| \leq 1$ as $k = 3$ and $N[y] = Y$. If $|N(x) \cap Y| \geq 2$ for some $x \in X$, then we take $D = \{x\}$. Suppose $|N(x) \cap Y| \leq 1$ for each $x \in X$. Since $G$ is connected, $N(x_1) \cap Y = \{y_1\}$ for some $x_1 \in X$ and $y_1 \in Y$. Let $x_2, x_3, y_2, y_3 \in V(G)$ such that $X = \{x_1, x_2, x_3\}$ and $Y = \{y_1, y_2, y_3\}$. If $G[v,x_2,x_3] \not \simeq K_3$, then we take $D = \{y_1\}$. Suppose $G[v,x_2,x_3] \simeq K_3$. If $x_iy_j \notin E(G)$ for every $i, j \in \{2, 3\}$, then we take $D = \{x_1\}$. If $x_iy_j \in E(G)$ for some $i, j \in \{2, 3\}$, then, since $\{v, x_2, x_3, y_j\} \subseteq N[x_i]$ and $N[x_1] \cap Y = \{y_1\}$, we take $D = \{x_i\}$.~\hfill{$\Box$}
\\  

\noindent
\textbf{Proof of Theorem~\ref{result k=3}.} The equality $\iota_3(B_{n,K_3}) = \left \lfloor \frac{n}{4} \right \rfloor$ is obtained similarly to the equality $\iota_2(B_{n,P_3}') = \left \lfloor\frac{4n-r}{14} \right \rfloor$ in Theorem~\ref{result k=2}.

We now prove the bound in the theorem, using the inductive argument in the proof of Theorem~\ref{result k=2}. Since $\iota_2(G)$ is an integer, it suffices to prove that $\iota_3(G) \leq n/4$. If $n \leq 7$, then the result is given by Lemma~\ref{small n k=3}. Suppose $n \geq 8$. Let $k = \max \{d(v) \colon v \in V(G)\}$. Since $G$ is connected, $k \geq 2$. If $k = 2$, then $G \simeq P_n$ or $G \simeq C_n$. Clearly, $\{5i \colon 1 \leq i \leq n/5\}$ is an $\mathcal{E}_3$-isolating set of $P_n$, and $\{6i-5 \colon 1 \leq i \leq \lfloor (n+5)/6 \rfloor\}$ is an $\mathcal{E}_3$-isolating set of $C_n$. Consequently, $\iota_3(G) \leq n/4$ if $k = 2$. Suppose $k \geq 3$. Let $v \in V(G)$ such that $d(v) = k$. Then, $|N[v]| \geq 4$. If $V(G) = N[v]$, then $\{v\}$ is an $\mathcal{E}_3$-isolating set of $G$, so $\iota_{3}(G) = 1 \leq n/4$. Suppose $V(G) \neq N[v]$. Let $G' = G-N[v]$ and $n' = |V(G')|$. Then, $n \geq n' + 4$ and $V(G') \neq \emptyset$. Let $\mathcal{H}$ be the set of components of $G'$. Let 
\begin{equation} \mathcal{H}' = \{H \in \mathcal{H} \colon \iota_3(H) > |V(H)|/4\}. \nonumber 
\end{equation} \vskip 0.1in

Suppose $\mathcal{H}' = \emptyset$. By Lemma~\ref{lemma} (with $X = \{v\}$ and $Y = N[v]$) and Lemma~\ref{lemmacomp},
\begin{align}
\iota_3(G) &\leq 1 + \iota_3(G') = 1 + \sum_{H \in \mathcal{H}} \iota_3(H) \leq \frac{|N[v]|}{4} + \sum_{H \in \mathcal{H}} \frac{|V(H)|}{4} = \frac{n}{4}. \nonumber 
\end{align} \vskip 0.1in

Now suppose $\mathcal{H}' \neq \emptyset$. By the induction hypothesis, each member of $\mathcal{H}'$ is a copy of $K_3$ or of $C_7$. Let $\mathcal{H}^{(1)} = \{H \in \mathcal{H}' \colon H \simeq K_3\}$ and $\mathcal{H}^{(2)} = \{H \in \mathcal{H}' \colon H \simeq C_7\}$. Thus, $\mathcal{H}' = \mathcal{H}^{(1)} \cup \mathcal{H}^{(2)}$.\medskip

For any $H \in \mathcal{H}$ and any $x \in N(v)$ such that $xy_{x,H} \in E(G)$ for some $y_{x,H} \in V(H)$, we say that $H$ is \emph{linked to $x$} and that $x$ is \emph{linked to $H$}. Since $G$ is connected, each member of $\mathcal{H}$ is linked to at least one member of $N(v)$. For each $x \in N(v)$, let $\mathcal{H}_x = \{H \in \mathcal{H} \colon H \mbox{ is linked to } x\}$, $\mathcal{H}'_x = \{H \in \mathcal{H}' \colon H \mbox{ is linked to } x\}$, and $\mathcal{H}_x^* = \{H \in \mathcal{H} \backslash \mathcal{H}' \colon H \mbox{ is linked to $x$ only}\}$. For each $H \in \mathcal{H} \backslash \mathcal{H}'$, let $D_H$ be an $\mathcal{E}_3$-isolating set of $H$ of size $\iota_3(H)$.\medskip 

For any $x \in N(v)$ and any $H \in \mathcal{H}_x'$, let $H_x' = H - \{y_{x,H}\}$. We have $H \simeq K_3$ or $H \simeq C_7$. If $H \simeq K_3$, then we take $D_{x,H} = \emptyset$. If $H \simeq C_7$, then $H$ is a $7$-cycle $(\{y_1, \dots, y_7\}, \{y_1y_2, \dots, y_6y_7, y_7y_1\})$ such that $y_{x,H} = y_1$, and we take $D_{x,H} = \{y_4\}$. Clearly, $D_{x,H}$ is an $\mathcal{E}_3$-isolating set of $H_x'$, and $|D_{x,H}| \leq (|V(H_x')| - 2)/4$. 
\\

\noindent
\emph{Case 1: $|\mathcal{H}'_x| \geq 2$ for some $x \in N(v)$.} For each $H \in \mathcal{H}' \backslash \mathcal{H}'_x$, let $x_H \in N(v)$ such that $H$ is linked to $x_H$. Let $X = \{x_H \colon H \in \mathcal{H}' \backslash \mathcal{H}'_x\}$. Note that $x \notin X$. Let 
\[D = \{v, x\} \cup X \cup \left( \bigcup_{H \in \mathcal{H}_x'} D_{x,H} \right) \cup \left( \bigcup_{H \in \mathcal{H}' \backslash \mathcal{H}_x'} D_{x_H,H} \right) \cup \left( \bigcup_{H \in \mathcal{H} \backslash \mathcal{H}'} D_{H} \right).\]
We have $V(G) = N[v] \cup \bigcup_{H \in \mathcal{H}} V(H)$, $y_{x,H} \in N[x]$ for each $H \in \mathcal{H}'_x$, and $y_{x_H,H} \in N[x_H]$ for each $H \in \mathcal{H}' \backslash \mathcal{H}'_x$, so $D$ is an $\mathcal{E}_3$-isolating set of $G$. Let 
\[Y = \{v, x\} \cup X \cup \{y_{x,H} \colon H \in \mathcal{H}_x'\} \cup \{y_{x_H,H} \colon H \in \mathcal{H}' \backslash \mathcal{H}'_x\}.\]
Since $|\mathcal{H}'_x| \geq 2$ and $|X| \leq |\mathcal{H}' \backslash \mathcal{H}_x'|$, $|Y| \geq 4 + 2|X| \geq 8 + 4|X| - 2|\mathcal{H}'_x| - 2|\mathcal{H}' \backslash \mathcal{H}'_x|$. We have
\begin{align} \iota_3(G) &\leq |D| = 2 + |X| + \sum_{H \in \mathcal{H}_x'} |D_{x,H}| + \sum_{H \in \mathcal{H}' \backslash \mathcal{H}_x'} |D_{x_H,H}| + \sum_{H \in \mathcal{H} \backslash \mathcal{H}'} |D_{H}| \nonumber \\
&\leq \frac{8 + 4|X|}{4} + \sum_{H \in \mathcal{H}_x'} \frac{|V(H_x')| - 2}{4} + \sum_{H \in \mathcal{H}' \backslash \mathcal{H}_{x}'} \frac{|V(H_{x_H}')| - 2}{4} + \sum_{H \in \mathcal{H} \backslash \mathcal{H}'} \frac{|V(H)|}{4} \nonumber \\
&\leq \frac{|Y|}{4} + \sum_{H \in \mathcal{H}_x'} \frac{|V(H_x')|}{4} + \sum_{H \in \mathcal{H}' \backslash \mathcal{H}_{x}'} \frac{|V(H_{x_H}')|}{4} + \sum_{H \in \mathcal{H} \backslash \mathcal{H}'} \frac{|V(H)|}{4} \leq \frac{n}{4}. \nonumber
\end{align}  \vskip 0.2in

\noindent
\emph{Case 2:}
\begin{equation} |\mathcal{H}'_x| \leq 1 \mbox{ \emph{for each} } x \in N(v). \label{k=3Hx<2} 
\end{equation} 
For each $H \in \mathcal{H}'$, let $x_H \in N(v)$ such that $H$ is linked to $x_H$. Let $X = \{x_H \colon H \in \mathcal{H}'\}$. Then, $\mathcal{H}' = \bigcup_{x \in X} \mathcal{H}_x'$. By (\ref{k=3Hx<2}), no two members of $\mathcal{H}'$ are linked to the same neighbour of $v$, so $|X| = |\mathcal{H}'|$. 
Let $W = N(v) \backslash X$, $W' = N[v] \backslash X$, and $Y_X = X \cup \{y_{x_H,H} \colon H \in \mathcal{H}'\}$. We have $|Y_X| = |X| + |\mathcal{H}'| = 2|X|$. Let 
\[D = \{v\} \cup X \cup \left( \bigcup_{H \in \mathcal{H}'}D_{x_H,H} \right) \cup \left( \bigcup_{H \in \mathcal{H} \backslash \mathcal{H}'} D_{H} \right).\] 
We have $V(G) = N[v] \cup \bigcup_{H \in \mathcal{H}} V(H)$ and $y_{x_H,H} \in N[x_H]$ for each $H \in \mathcal{H}$, so $D$ is an $\mathcal{E}_3$-isolating set of $G$.\medskip

\noindent
\emph{Subcase 2.1: $|W| \geq 3$.} Similarly to Case 1, we have
\begin{align} \iota_3(G) &\leq |D| = 1 + |X| + \sum_{H \in \mathcal{H}'} |D_{x_H,H}| + \sum_{H \in \mathcal{H} \backslash \mathcal{H}'} |D_{H}| \nonumber \\
&\leq \frac{|\{v\} \cup W| + 4|X|}{4} + \sum_{H \in \mathcal{H}'} \frac{|V(H_{x_H}')| - 2}{4} + \sum_{H \in \mathcal{H} \backslash \mathcal{H}'} \frac{|V(H)|}{4} \nonumber \\
&= \frac{|W'| + |Y_X|}{4} + \sum_{H \in \mathcal{H}'} \frac{|V(H_{x_H}')|}{4} + \sum_{H \in \mathcal{H} \backslash \mathcal{H}'} \frac{|V(H)|}{4} = \frac{n}{4}. \nonumber   
\end{align} \vskip 0.1in

\noindent
\emph{Subcase 2.2: $|W| \leq 2$.}\medskip

\noindent
\emph{Subsubcase 2.2.1: Some member $H$ of $\mathcal{H}'$ is linked to $x_H$ only.} Let $G^* = G - (\{x_H\} \cup V(H))$. Then, $G^*$ has a component $G_v^*$ such that $N[v] \backslash \{x_H\} \subseteq V(G_v^*)$, and the other components of $G^*$ are the members of $\mathcal{H}_{x_H}^*$. Let $D^*$ be an $\mathcal{E}_3$-isolating set of $G_v^*$ of size $\iota_3(G_v^*)$. Let $D' = D^* \cup \{x_H\} \cup D_{x_H,H} \cup \bigcup_{I \in \mathcal{H}_{x_H}^*} D_I$. Since $D'$ is an $\mathcal{E}_3$-isolating set of $G$,
\begin{equation} \iota_3(G) \leq |D^*| + 1 + |D_{x_H,H}| + \sum_{I \in \mathcal{H}_{x_H}^*} |D_I| \leq \iota_3(G_v^*) + \frac{|\{x_H\} \cup V(H)|}{4} + \sum_{I \in \mathcal{H}_{x_H}^*} \frac{|V(I)|}{4}.  \nonumber
\end{equation}
This yields $\iota_3(G) \leq n/4$ if $\iota_3(G_v^*) \leq |V(G_v^*)|/4$. Suppose $\iota_3(G_v^*) > |V(G_v^*)|/4$. By the induction hypothesis, $G_v^* \simeq K_3$ or $G_v^* \simeq C_7$. We have $v \in N[x_H]$. Thus, if $G_v^* \simeq K_3$, then $D' \backslash D^*$ is an $\mathcal{E}_3$-isolating set of $G$, and hence
\[\iota_3(G) \leq \frac{|\{x_H\} \cup V(H)|}{4} + \sum_{I \in \mathcal{H}_{x_H}^*} \frac{|V(I)|}{4} < \frac{n}{4}.\] 
Suppose $G_v^* \not\simeq K_3$. Then, $G_v^*$ is a $7$-cycle $(\{y_1, \dots, y_7\}, \{y_1y_2, \dots, y_6y_7, y_7y_1\})$ such that $v = y_1$. Thus, since $v \in N[x_H]$, $(D' \backslash D^*) \cup \{y_4\}$ is an $\mathcal{E}_3$-isolating set of $G$, and hence 
\begin{equation} \iota_3(G) \leq 1 + \frac{|\{x_H\} \cup V(H)|}{4} + \sum_{I \in \mathcal{H}_{x_H}^*} \frac{|V(I)|}{4} < \frac{|V(G_v^*)|}{4} + \frac{n - |V(G_v^*)|}{4} = \frac{n}{4}. \nonumber
\end{equation}  \vskip 0.1in

\noindent
\emph{Subsubcase 2.2.2: For each $H \in \mathcal{H}'$, $H$ is linked to some $x_H' \in N(v) \backslash \{x_H\}$.} Let $X' = \{x_H' \colon H \in \mathcal{H}'\}$. By (\ref{k=3Hx<2}), for every $H, I \in \mathcal{H}'$ with $H \neq I$, we have $x_H' \neq x_I'$ and $x_H \neq x_I'$. Thus, $|X'| = |\mathcal{H}'|$ and $X' \subseteq W$. This yields $|\mathcal{H}'| \leq 2$ as $|W| \leq 2$.\medskip

Suppose $|\mathcal{H}'| = 2$. Then, $|W| = 2$ and $X' = W$. Let $H_1$ and $H_2$ be the two members of $\mathcal{H}'$. Each of $H_1$ and $H_2$ is a $\{K_3,C_7\}$-graph. Let $x_1 = x_{H_1}$, $x_2 = x_{H_2}$, $x_1' = x_{H_1}'$, and $x_2' = x_{H_2}'$. We have $N(v) = \{x_1, x_2, x_1', x_2'\}$. Let $Y_1 = N[y_{x_1, H_1}] \cap (\{x_1\} \cup V(H_1))$. Then, $x_1 \in Y_1$ and $|Y_1| = 4$. Let $G^* = G - Y_1$. Then, $G^*$ has a component $G_v^*$ such that $N[v] \backslash \{x_1\} \subseteq V(G_v^*)$, and the members of $\mathcal{H}_{x_1}^*$ are also components of $G^*$. If $G_v^*$ and the members of $\mathcal{H}_{x_1}^*$ are the only components of $G^*$, then we take $D_1 = \{y_{x_1, H_1}\}$; otherwise, $H_1 - N_{H_1}[y_{x_1, H_1}]$ is the only other component of $G^*$, $H_1 \simeq C_7$, $H_1 - N_{H_1}[y_{x_1, H_1}]\simeq P_4$, and we take $D_1 = \{y_{x_1, H_1}\} \cup D_{x_1, H_1}$. Let $D^*$ be an $\mathcal{E}_3$-isolating set of $G_v^*$ of size $\iota_3(G_v^*)$. Then, $D^* \cup D_1 \cup \bigcup_{H \in \mathcal{H}_{x_1}^*} D_H$ is an $\mathcal{E}_3$-isolating set of $G$. Since $d_{G_v^*}(v) = 3$, $G_v^*$ is not a $\{K_3, C_7\}$-graph. By the induction hypothesis, $\iota_3(G_v^*) \leq |V(G_v^*)|/4$. We have
\begin{align} \iota_3(G) &\leq |D^*| + |D_1| + \sum_{H \in \mathcal{H}_{x_1}^*} |D_H| \nonumber \\
&\leq \frac{|V(G_v^*)|}{4} + \frac{|(\{x_1\} \cup V(H_1)) \backslash V(G_v^*)|}{4} + \sum_{H \in \mathcal{H}_{x_1}^*} \frac{|V(H)|}{4} = \frac{n}{4}.   \nonumber
\end{align}  \vskip 0.1in

Finally, suppose $|\mathcal{H}'| = 1$. Let $H_1$ be the member of $\mathcal{H}'$, and let $x_1 = x_{H_1}$ and $x_1' = x_{H_1}'$. We have $3 \leq d(v) = |X| + |W| = 1 + |W| \leq 3$, so $d(v) = 3$, and hence $N(v) = \{x_1, x_1', w\}$ for some $w \in V(G) \backslash \{x_1, x_1'\}$. Let $Y = N[v] \cup V(H_1)$. Then, 
\begin{equation} V(G) = V(G[Y]) \cup \bigcup_{H \in \mathcal{H} \backslash \mathcal{H}'} V(H). \label{k=3GY}
\end{equation}

Suppose $H_1 \simeq C_7$. We first show that $G[Y]$ has an $\mathcal{E}_3$-isolating set $D_Y$ of size~$2$. We may assume that $H_1 = C_7$. For each $i \in [7]$, let $y_i$ be the vertex $i$ of $C_7$. We may assume that $y_{x_1,H_1} = y_1$. Let $Y_1 = \{x_1, y_1, y_2, y_7\}$. We have $Y_1 \subseteq N[y_1]$. Let $I = G[Y] - Y_1$. Then, $V(I) = \{v, x_1', w, y_3, y_4, y_5, y_6\}$ and $vx_1', vw, y_3y_4, y_4y_5, y_5y_6 \in E(I)$. If $I$ is connected and $I \not\simeq C_7$, then, by Lemma~\ref{small n k=3}, $I$ has an $\mathcal{E}_3$-isolating set $D_{I}$ of size $1$, so we take $D_Y = \{y_1\} \cup D_{I}$. Suppose that $I$ is not connected. Then, $E(I) \backslash \{x_1'w\} = \{vx_1', vw, y_3y_4, y_4y_5, y_5y_6\}$, and hence $N(y_i) \cap \{x_1', w\} = \emptyset$ for each $i \in \{3, 4, 5, 6\}$. If $N(y) \cap \{x_1', w\} \neq \emptyset$ for some $y \in \{y_2, y_7\}$, then we take $D_Y = \{y_1, y\}$. If $N(y_2) \cap \{x_1', w\} = N(y_7) \cap \{x_1', w\} = \emptyset$, then we take $D_Y = \{x_1, y_4\}$. Now suppose $I \simeq C_7$. We may assume that $E(I) = \{vx_1', x_1'y_3, y_3y_4, y_4y_5, y_5y_6, y_6w, wv\}$. Then, we take $D_Y = \{y_3, y_6\}$ as $G[Y] - N[\{y_3, y_6\}] = \{v, x_1, y_1\}$ and $vy_1 \notin E(G)$.\medskip

Therefore, $\iota_2(G[Y]) = 2$. If $V(G) = Y$, then $G = G[Y]$, so $\iota_3(G) = 2 < n/4$. Suppose $V(G) \neq Y$. Then, $\mathcal{H} \backslash \mathcal{H}' \neq \emptyset$. Thus, $uz \in E(G)$ for some $u \in N(v)$ and some $z \in V(G) \backslash Y$. 
Let $u' \in \{x_1, x_1'\} \backslash \{u\}$. Suppose $u' = x_1$.  Suppose that $I$ is connected. Let $G^* = G - Y_1$. Then, $G^*$ has a component $G_v^*$ such that $V(I) \cup \{z\} \subseteq V(G_v^*)$, and the other components of $G^*$ are the members of $\mathcal{H}_{x_1}^*$. We have $|V(G_v^*)| \geq 8$. By the induction hypothesis, $\iota_3(G_v^*) \leq |V(G_v^*)|/4$. Let $D^*$ be an $\mathcal{E}_3$-isolating set of $G_v^*$ of size $\iota_3(G_v^*)$. Then, $D^* \cup \{y_1\} \cup \bigcup_{H \in \mathcal{H}_{x_1}^*} D_H$ is an $\mathcal{E}_3$-isolating set of $G$, and hence
\begin{equation} \iota_3(G) \leq |D^*| + 1 + \sum_{H \in \mathcal{H}_{x_1}^*} |D_H| \leq \frac{|V(G_v^*)|}{4} + \frac{|Y_1|}{4} + \sum_{H \in \mathcal{H}_{x_1}^*} \frac{|V(H)|}{4} = \frac{n}{4}.  \nonumber
\end{equation}
Now suppose that $I$ is not connected. Then, $xy_i \notin E(G)$ for each $x \in N(v)$ and each $i \in \{3, 4, 5, 6\}$. Thus, $y_{x_1',H_1} \in \{y_1, y_2, y_7\}$. Since $Y_1 \subseteq N[y_1]$, $d(y_1) \geq 3$. Since $k = d(v) = 3$, we obtain $d(y_1) = 3$, so $N[y_1] = Y_1$. Thus, $x_1'y_1 \notin E(G)$, and hence $y_{x_1',H_1} \in \{y_2, y_7\}$. Let $G^* = G - (Y \backslash \{v, w\})$. Then, $G^*$ has a component $G_v^*$ such that $\{v, w\} \subseteq V(G_v^*)$, and the other components of $G^*$ are the members of $\mathcal{H} \backslash (\mathcal{H}' \cup \mathcal{H}_w)$. If $G_v^*$ is a $\{K_3, C_7\}$-graph, then $\iota_3(G_v^*) = \frac{|V(G_v^*)| + 1}{4}$; otherwise, $\iota_3(G_v^*) \leq \frac{|V(G_v^*)|}{4}$ by the induction hypothesis. Let $D^*$ be an $\mathcal{E}_3$-isolating set of $G_v^*$ of size $\iota_3(G_v^*)$. Then, $D^* \cup \{y_1, y_{x_1',H_1}\} \cup \bigcup_{H \in \mathcal{H} \backslash (\mathcal{H}' \cup \mathcal{H}_w)} D_H$ is an $\mathcal{E}_3$-isolating set of $G$, and hence
\begin{equation} \iota_3(G) \leq \frac{|V(G_v^*)| + 1}{4} + \frac{|Y \backslash \{v, w\}| - 1}{4} + \sum_{H \in \mathcal{H} \backslash (\mathcal{H}' \cup \mathcal{H}_w)} \frac{|V(H)|}{4} = \frac{n}{4}. \nonumber
\end{equation}
Similarly, $\iota_3(G) \leq n/4$ if $u' = x_1'$.\medskip

Now suppose $H_1 \not\simeq C_7$. Then, $H_1 \simeq K_3$. We have $Y = \{v, x_1, x_1', w, y_1, y_2, y_3\}$, where $\{y_1, y_2, y_3\} = V(H_1)$ and $x_1y_1 \in E(G)$. Let $Y_1 = \{x_1\} \cup V(H_1)$. Similarly to the above, since $k = d(v) = 3$, we have $N[y_1] = Y_1$ and $y_{x_1',H_1} \in \{y_2, y_3\}$. Since $n \geq 8$, $uz \in E(G)$ for some $u \in N(v)$ and some $z \in V(G) \backslash Y$. Let $u' \in \{x_1, x_1'\} \backslash \{u\}$. As above, we may assume that $u' = x_1$. Let $G^* = G - Y_1$. Then, $G^*$ has a component $G_v^*$ such that $\{v, x_1', w, z\} \subseteq V(G_v^*)$, and the other components of $G^*$ are the members of $\mathcal{H}_{x_1}^*$. Suppose $G_v^* \not\simeq C_7$. By the induction hypothesis, $\iota_3(G_v^*) \leq |V(G_v^*)|/4$. Let $D^*$ be an $\mathcal{E}_3$-isolating set of $G_v^*$ of size $\iota_3(G_v^*)$. Then, $D^* \cup \{y_1\} \cup \bigcup_{H \in \mathcal{H}_{x_1}^*} D_H$ is an $\mathcal{E}_3$-isolating set of $G$, and hence
\begin{equation} \iota_3(G) \leq |D^*| + 1 + \sum_{H \in \mathcal{H}_{x_1}^*} |D_H| \leq \frac{|V(G_v^*)|}{4} + \frac{|Y_1|}{4} + \sum_{H \in \mathcal{H}_{x_1}^*} \frac{|V(H)|}{4} = \frac{n}{4}.  \nonumber
\end{equation}
Now suppose $G_v^* \simeq C_7$. Then, for some distinct $z_1, z_2, z_3, z_4 \in V(G) \backslash Y$, $V(G_v^*) = \{v, x_1', w, z_1, z_2, z_3, z_4\}$ and $E(G_v^*) = \{vx_1', x_1'z_1, z_1z_2, z_2z_3, z_3z_4, z_4w, wv\}$. Since $y_{x_1',H_1} \in N(x_1')$, we have that $\{v, x_1'\} \cup \bigcup_{H \in \mathcal{H}_{x_1}^*} D_H$ is an $\mathcal{E}_3$-isolating set of $G$, so
$\iota_3(G) < (|Y_1| + |V(G_v^*)|)/4 + \sum_{H \in \mathcal{H}_{x_1}^*} |V(H)|/4 = n/4$. \hfill{$\Box$}

\section{Problems}
We conclude this paper by posing some problems and a conjecture.

We may assume that the vertex set of an $n$-vertex graph is $[n]$. Let 
\[\mathcal{G} = \{G \colon G \mbox{ is a connected graph, } V(G) = [n] \mbox{ for some } n \geq 1\}.\] 
%
%
%
Thus, $\mathcal{G}$ is an infinite set. For any integer $k \geq 1$ and any real number $\alpha > 0$, let 
\[\mathcal{G}(k,\alpha) = \{G \in \mathcal{G} \colon \iota_k(G) \leq \lfloor \alpha |V(G)| \rfloor\},\] 
and let
\[\mathcal{G}(k,\alpha)^* = \{G \in \mathcal{G}(k,\alpha) \colon \iota_k(G) = \lfloor \alpha |V(G)| \rfloor\} \quad \mbox{and} \quad \mathcal{G}[k,\alpha] = \mathcal{G} \backslash \mathcal{G}(k,\alpha).\]
Thus, $\mathcal{G}[k,\alpha] = \{G \in \mathcal{G} \colon \iota_k(G) > \lfloor \alpha |V(G)| \rfloor\}$. In view of Theorems~\ref{CHresult}--\ref{result k=3}, we pose the following problem.

\begin{problem} \label{problem1} (a) Given $k \geq 1$, is there a rational number $c_k$ such that $\mathcal{G}{[k,c_k]}$ is finite and $\mathcal{G}{(k,c_k)}^*$ is infinite? \\
(b) If $c_k$ exists, then determine $c_k$, $\mathcal{G}{[k,c_k]}$, and (at least) an infinite subset of $\mathcal{G}(k,c_k)^*$.
\end{problem}
By Theorem~\ref{CHresult}, $c_1 = 1/3$, $\mathcal{G}{[1,c_1]} = \{K_2\} \cup \{G \in \mathcal{G} \colon G \simeq C_5\}$, and $\{G \in \mathcal{G} \colon G \simeq B_{n,K_2} \mbox{ for some } n \geq 3\} \subseteq \mathcal{G}{(1,c_1)^*}$. By Theorem~\ref{result k=3}, $c_3 = 1/4$, $\mathcal{G}[3,c_3] = \{K_3\} \cup \{G \in \mathcal{G} \colon G \simeq C_7\}$, and $\{G \in \mathcal{G} \colon G \simeq B_{n,K_3} \mbox{ for some } n \geq 1\} \subseteq \mathcal{G}(3,c_3)^*$. Now consider $k = 2$. By Theorem~\ref{result k=2}, $\iota_2(G) \leq \lfloor 4|V(G)|/14 \rfloor = \lfloor 2|V(G)|/7 \rfloor$ for any $G \in \mathcal{G}$ such that $G$ is not a $\{P_3, K_3, C_6\}$-graph. 
%
%
%
For $r \geq 1$, let $F_1', \dots, F_r'$ be pairwise vertex-disjoint copies of $C_6'$ such that, for each $i \in [r]$, $i$ is the leaf of $F_i'$, and let $B_{7r,C_6}'$ be the union of $F_1', \dots, F_r'$, and $P_r$. Thus, $B_{7r,C_6}'$ is a connected $7r$-vertex graph. For each $i \in [r]$, let $v_i$ be the neighbour of $i$ in $F_i'$, and let $v_i' \in V(F_i') \backslash N_{F_i'}(v_i)$. Clearly, $\bigcup_{i = 1}^{r} \{i,v_i'\}$ is an $\mathcal{E}_2$-isolating set of $B_{7r,C_6}'$, and for every $\mathcal{E}_2$-isolating set $D$ of $B_{7r,C_6}'$, $|D \cap V(F_i')| \geq 2$ for each $i \in [r]$. Thus, $\iota_2(B_{7r,C_6}') = 2r = 2|V(B_{7r,C_6}')|/7$. Therefore, $c_2 = 2/7$, $\mathcal{G}{[2,c_2]} = \{G \in \mathcal{G} \colon G  \mbox{ is a } \{P_3, K_3, C_6\}\mbox{-graph}\}$, and $\{G \in \mathcal{G} \colon G \simeq B_{7r,C_6}' \mbox{ for some } r \geq 1\} \subseteq \mathcal{G}{(2,c_2)^*}$.

\begin{conj} For $k \geq 1$, $c_k$ exists.
\end{conj}

Our next problem is stronger than Problem~\ref{problem1}(a).

\begin{problem} Given $k$ and $n$, what is the smallest rational number $c_{k,n}$ such that $\iota_k(G) \leq c_{k,n} n$ for every $n$-vertex graph $G$?
\end{problem}
By Theorem~\ref{CHresult}, $c_{1,2} = 1/2$, $c_{1,5} = 2/5$, and $c_{1,n} = \lfloor n/3 \rfloor / n$ for any $n \in \mathbb{N} \backslash \{2, 5\}$.  By Theorem~\ref{result k=2} and the above,  $c_{2,3} = c_{2,6} = 1/3$ and $c_{2,n} = \lfloor 2n/7 \rfloor / n$ for any $n \in \mathbb{N} \backslash \{3, 6\}$. By Theorem~\ref{result k=3}, $c_{3,3} = 1/3$, $c_{3,7} = 2/7$, and $c_{3,n} = \lfloor n/4 \rfloor / n$ for any $n \in \mathbb{N} \backslash \{3, 7\}$.


\begin{thebibliography}{}

\bibitem{Borg1} P. Borg, Isolation of cycles, Graphs and Combinatorics 36 (2020), 631--637.

\bibitem{BFK} P. Borg, K. Fenech and P. Kaemawichanurat, Isolation of $k$-cliques, Discrete Mathematics 343 (2020), paper 111879.

\bibitem{CaHa17} Y. Caro and A. Hansberg, Partial Domination - the Isolation Number of a Graph, Filomat 31:12 (2017), 3925--3944. 




\bibitem{C} E.J. Cockayne, Domination of undirected graphs -- A survey, Lecture Notes in Mathematics, Volume 642, Springer, 1978, 141--147.

\bibitem{CH} E.J. Cockayne and S.T. Hedetniemi, Towards a theory of domination in graphs, Networks 7 (1977), 247--261.



\bibitem{HHS} T.W. Haynes, S.T. Hedetniemi and P.J. Slater, Fundamentals of Domination in Graphs, Marcel Dekker, Inc., New York, 1998.

\bibitem{HHS2} T.W. Haynes, S.T. Hedetniemi and P.J. Slater (Editors), Domination in Graphs: Advanced Topics, Marcel Dekker, Inc., New York, 1998.

\bibitem{HL} S.T. Hedetniemi and R.C. Laskar (Editors), Topics on Domination, Discrete Mathematics, Volume 86, 1990.

\bibitem{HL2} S.T. Hedetniemi and R.C. Laskar, Bibliography on domination in graphs and some basic definitions of domination parameters, Discrete Mathematics 86 (1990), 257--277.




\bibitem{Ore} O. Ore, Theory of graphs, American Mathematical Society Colloquium Publications, Volume 38, American Mathematical Society, Providence, R.I., 1962.



\bibitem{West} D.B. West, Introduction to graph theory (second edition), Prentice Hall, 2001.
	
\end{thebibliography}
\end{document}